\documentclass [11pt] {article}

\usepackage {amssymb}

\title {Homogeneous Fedosov Star Products\\
        on Cotangent Bundles I: \\
        Weyl and Standard Ordering with Differential Operator
        Representation}

\author {{\bf
          Martin
          Bordemann\thanks{Martin.Bordemann@physik.uni-freiburg.de}~,
          \addtocounter{footnote}{2}
          Nikolai
          Neumaier\thanks{Nikolai.Neumaier@physik.uni-freiburg.de}~,
          Stefan
          Waldmann\thanks{Stefan.Waldmann@physik.uni-freiburg.de}
         } \\[3mm]
         Fakult\"at f\"ur Physik\\Universit\"at Freiburg \\
         Hermann-Herder-Str. 3 \\
         79104 Freiburg i.~Br., F.~R.~G \\[3mm]
        }

\date{FR-THEP-97/10 \\[1mm]
      July 1997\\[1mm]
     }

\newcommand {\W}   {\mathcal W}
\newcommand {\WL}  {\mbox{$\mathcal W \! \otimes \! \Lambda$}}
\newcommand {\WLX} {\mbox{$\mathcal W \! \otimes \! \Lambda
                          \otimes \!\mathcal X$}}

\newcommand {\degs} {{\rm deg}_s}
\newcommand {\dega} {{\rm deg}_a}
\newcommand {\degh} {{\rm deg}_\lambda}
\newcommand {\Deg} {{\rm Deg}}

\newcommand {\TinyW} {{\mbox{\rm \tiny W}}}

\newcommand {\starw} {*_{\mbox{\rm \tiny W}}}

\newcommand {\wrep} {\varrho_\TinyW}

\newcommand {\TinyF} {{\mbox{\rm \tiny F}}}
\newcommand {\fed} {\circ_\TinyF}
\newcommand {\adfed} {{\rm ad}_\TinyF}
\newcommand {\rfed} {r_\TinyF}
\newcommand {\Dfed} {\mathcal D_\TinyF}
\newcommand {\taufed} {\tau_\TinyF}
\newcommand {\starf} {*_{\mbox{\rm \tiny F}}}
\newcommand {\Dfedcl} {\mathcal D_\TinyF^{\rm cl}}
\newcommand {\taufedcl} {\tau_\TinyF^{\rm cl}}
\newcommand {\fibfrep} {\tilde{\varrho}_\TinyF}

\newcommand {\TinyS} {{\mbox{\rm \tiny S}}}
\newcommand {\std} {\circ_\TinyS}
\newcommand {\ads} {{\rm ad}_\TinyS}
\newcommand {\rstd} {r_\TinyS}
\newcommand {\Dstd} {\mathcal D_\TinyS}
\newcommand {\taustd} {\tau_\TinyS}
\newcommand {\stars} {*_{\mbox{\rm \tiny S}}}
\newcommand {\fibsrep} {\tilde{\varrho}_\TinyS}
\newcommand {\srep} {\varrho_\TinyS}

\newcommand {\CCs} {{\mathcal C_\TinyS}}
\newcommand {\Cs}  {{C_\TinyS}}
\newcommand {\Deltafib} {{\Delta_{\mbox{\rm \tiny fib}}}}

\newcommand {\starg} {*_{\mbox{\rm \tiny G}}}

\newcommand {\BEQ} [1] {\begin {equation} \label {#1}}
\newcommand {\EEQ} {\end {equation}}

\newcommand {\Lie} {\mathcal L}
\newcommand {\cc} [1] {\overline {{#1}}}

\newcommand {\id} {{\mathsf {id}}}
\newcommand {\tr} {{\mathsf {tr}}}
\newcommand {\ad} {{\mathrm {ad}}}
\newcommand {\End} {{\mathsf {End}}}
\newcommand {\Div} {{\mathsf {div}}}
\newcommand {\Ric} {{\mathrm {Ric}}}

\DeclareMathAlphabet{\mathbfit}{OT1}{cmr}{bx}{it}
\newcommand {\im} {{\mathbfit i}}

\newenvironment {PROOF}{\small {\sc Proof:}}{{\hspace*{\fill} $\square$}}

\newtheorem {LEMMA} {Lemma} [section]
\newtheorem {PROPOSITION} [LEMMA] {Proposition}
\newtheorem {THEOREM} [LEMMA] {Theorem}
\newtheorem {COROLLARY} [LEMMA] {Corollary}
\newtheorem {DEFINITION}[LEMMA] {Definition}

\begin {document}

\maketitle

\begin {abstract}

In this paper we construct homogeneous star products of Weyl type on 
every cotangent bundle $T^*Q$ by means of the Fedosov procedure using
a symplectic torsion-free connection on $T^*Q$ homogeneous of
degree zero with respect to the Liouville vector field. 
By a fibrewise equivalence transformation we construct a homogeneous
Fedosov star product of standard ordered type equivalent to the 
homogeneous Fedosov star product of Weyl type.
Representations for both star product algebras by differential operators on 
functions on $Q$ are constructed leading in the case of the standard
ordered product to the usual standard ordering prescription for
smooth complex-valued
functions on $T^*Q$ polynomial in the momenta (where an arbitrary fixed 
torsion-free connection $\nabla_0$ on $Q$ is used). 
Motivated by the flat case $T^*\mathbb R^n$ another homogeneous star 
product of Weyl type corresponding to the Weyl ordering prescription is 
constructed. The example of the cotangent bundle of an arbitrary Lie 
group is explicitly computed and the star product given by Gutt is
rederived in our approach.

\end {abstract}

\newpage


\section {Introduction}
\label {IntroSec}

The concept of deformation quantization defined in \cite{BFFLS78}
has now been well-established on every symplectic manifold:
existence of the formal associative deformation of the pointwise
multiplication of smooth functions, the so-called star product,
had been shown by \cite{DL83} and \cite{Fed94}, and their
classification up to equivalence transformations by formal power series
in the second de Rham cohomology group is due to \cite{NT95a,NT95b} and
\cite{BCG96}.

The symplectic manifolds which are mostly used by physicists are
cotangent bundles $T^*Q$ of a smooth manifold $Q$, the
configuration space of the classical dynamical system, which is rather
often taken to be an open subset of $\mathbb R^{2n}$. There is a large
amount of literature concerning star products on cotangent bundles
(cf.~e.~g.~\cite {CFS92}, \cite{DL83a}, \cite{Pfl95}), 
differential operators and their symbolic calculus 
(\cite{Und78},\cite{Wid78}), and also geometric quantization on 
cotangent bundles (see e.~g.~\cite{Woo80} and references therein).

The main motivation for us to write this paper was to apply the formal GNS
construction in deformation quantization developed by two of us
(cf.~\cite{BW96b}) to the particular case of $T^*Q$: this method
(which basically copies the standard GNS representation in the theory of
$C^*$-algebras) allows to construct formal pre-Hilbert space
representations of the associative algebra (over the field of formal
Laurent series with complex coefficients, $\mathbb C(\!(\lambda)\!)$) of
all formal Laurent series with coefficients in the space of all smooth
complex-valued functions on a symplectic manifold equipped with a star
product. The basic ingredient is a formally positive
$\mathbb C(\!(\lambda)\!)$-linear functional on this algebra.
In case $Q=\mathbb R^n$ the usual Weyl ordered Schr\"odinger
representation could thus be reconstructed by means of an integral
over configuration space (see \cite{BW96b}) as well as the ordinary
WKB expansion by means of a certain functional with support on a
projectable Lagrangean submanifold of $\mathbb R^{2n}$ contained in a
classical energy surface (see \cite{BW96c}).

When starting to work on general cotangent bundles we realized that we
had to develop first a good deal of compact practical formulas for
certain star products on $T^*Q$ and their possible representations
as formal series with coefficients in differential operators on $Q$
before we could start
checking that even the simplest functional which consists in integration
over $Q$ (with respect to some volume in case $Q$ is orientable) is
formally positive. Therefore on one hand this paper will simply prepare
the
grounds for a second paper (see \cite{BNW97b}) in which we shall define
formal GNS representations on $T^*Q$ and also compare our results with
those obtained by analytic techniques (see e.~g.~\cite {CFS92}, 
\cite {Pfl95}, \cite{Und78}, \cite{Wid78}). On the other hand
we feel that some of our results, viz. a fairly explicit Fedosov
construction on arbitrary $T^*Q$, and a rather simple closed formula
relating a star product based on standard ordering and a particular
star product of Weyl type (which is different from the `most natural'
Fedosov star product!) based on a generalization of Weyl ordering
may well be of independent interest and useful in computations
because the usual techniques of asymptotic expansions of certain
integrals in normal coordinates are not needed.

Before summarizing our results let us first motivate our programme
by the simple example $Q=\mathbb R^n$ in which everything can
explicitly be computed (see e.~g.~\cite{AW70}, \cite{BW96b}):

The quantization method which is frequently used by physicists for
$Q=\mathbb R^n$ in the Schr\"odinger picture proceeds as follows:
Except for relativistic phase space functions such as the energy
$\sqrt{m^2+p^2}$ of a free particle in $\mathbb R^n$ the great majority
of classical observables occurring in physics are
{\em polynomial in the momenta $p$}
i.~e.~smooth functions
$F:\mathbb R^{2n}\rightarrow \mathbb C:(q,p)\mapsto F(q,p)$ which take the
form\footnote{From now on we shall use the
Einstein summation convention where the sum over repeated coordinate
indices $i_r$ where mostly $1\leq i_r\leq n$ is automatic.}
\BEQ {polyflat}
    (q,p) \mapsto F(q,p) =
    \sum_{k=0}^N\frac{1}{k!} F_k^{i_1\cdots i_k}(q)
    p_{i_1} \cdots p_{i_k}
\EEQ
(where the $F_k^{i_1\cdots i_k}$ are smooth complex-valued functions on
$\mathbb R^n$). On the space of all these functions which we shall denote
by $C^\infty_{pp}(\mathbb R^{2n})$ (and which clearly forms a Poisson
subalgebra of $C^\infty(\mathbb R^{2n})$) one establishes a linear 
bijection to the space of all differential operators on 
functions $\psi$ on $\mathbb R^n$ according to the following rule: 
a smooth complex-valued function
$q \mapsto f(q)$ is mapped to the multiplication with $f$, i.~e. 
$\psi\mapsto f\psi$, the coordinates $p_l$ are mapped to $\frac{\hbar}{\im}
\frac{\partial}{\partial q^l}$ and for a general function polynomial in 
the momenta a so-called {\em ordering prescription} is applied
to extend the map to a bijection:
An important example is the {\em standard ordering prescription} where
a function of the above form (\ref{polyflat}) is mapped to its
{\em standard representation},
$\srep(F)$, in the following way:
\BEQ {StandOrdFlat}
    \srep(F)(\psi):q\mapsto \sum_{k=0}^N
    \frac{1}{k!} \left(\frac{\hbar}{\im}\right)^k
    F_k^{i_1\cdots i_k} (q)
    \frac{\partial^k\psi}{\partial q^{i_1}\cdots \partial q^{i_k}}(q) .
\EEQ
The standard representation, however, is unphysical in the sense that
the differential operators $\srep(F)$ are not symmetric
(when $F$ takes real values) on, say, the space
${\mathcal D}(\mathbb R^n)$ of all smooth complex-valued functions with
compact support with the standard inner product given by the
Lebesgue integral
(i.~e. $\langle \phi, \psi\rangle:=\int \cc{\phi(q)} \psi(q) d^n q$ where
$\cc{\rule{0mm}{0.5em} \;\;}$ denotes pointwise complex conjugation):
it is easily seen by induction and repeated partial integration that
the formal adjoint $\srep(F)^\dagger$ of $\srep(F)$ (i.~e.
$\langle \srep(F)^\dagger\phi,\psi\rangle=\langle\phi,\srep(F)\psi\rangle$)
is given by the differential operator
\BEQ {FormalAdjoint}
  \srep(F)^\dagger = \srep(N^2 \cc{F})
\EEQ
where
\BEQ {SFlat}
     N:=e^{\frac{\hbar}{2\im} \frac{\partial^2}
          {\partial q^k\partial p_k}}
\EEQ
is a well-defined bijective linear map on all the functions polynomial in
the momenta (\ref{polyflat}). These unphysical features can be remedied
by defining the {\em Weyl representation} of $F$ by
\BEQ {WeylundStandard}
   \wrep(F) := \srep(NF)
\EEQ
which is also a bijection and clearly gives
\BEQ {WeylAdjoint}
   \wrep(F)^\dagger = \wrep(\cc{F})
\EEQ
such that real-valued functions are now mapped to symmetric operators.
In case $F$ is a polynomial function in $p$ and in $q$ it is not hard
to see using the Baker-Campbell-Hausdorff series of the Heisenberg
Lie algebra spanned by $1,\srep(q^1), \ldots, \srep(q^n)$,
$\srep(p_1), \ldots, \srep(p_n)$ that $\wrep(F)$ can be obtained by
the so-called {\em Weyl ordering prescription} by means of which
monomials are mapped to totally symmetrized operators (see
e.~g.~\cite{AW70}, \cite{BW96b}), i.~e.
\BEQ {WeylOrdering}
    q^{i_1}\cdots q^{i_a}p_{j_1}\cdots p_{j_b}
       \mapsto \frac{1}{(a+b)!}\sum_{\sigma\in S_{a+b}}
                      A_{\sigma(1)}\cdots A_{\sigma(a+b)},
\EEQ
where $A_r:=\srep(q^{i_r})$ for all $1\leq r\leq a$ and
$A_r:=\srep(p_{j_{r-a}})$ for all $a+1\leq r\leq a+b$.

The usual Moyal-Weyl star product $\starw$ in $\mathbb R^{2n}$ of two
smooth complex-valued functions $F,G$ polynomial in the momenta,
\BEQ {Weylflat}
       (F\starw G) (q,p):=
       e^{\frac{{\im}\hbar}{2}
          \left(\frac{\partial^2}{\partial q^i\partial p'_i}
              - \frac{\partial^2}{\partial q'^i\partial p_i} \right)}
       F(q,p)G(q',p')\bigg|_{q=q',p=p'}
\EEQ
can for instance be obtained from the multiplication of the two Weyl
representations, i.~e.:
\BEQ {MoyalWeylProd}
          \wrep(F\starw G)=\wrep(F)\wrep(G).
\EEQ
Likewise, the multiplication of the two standard representations
of $F$ and $G$ gives rise to another
{\em star product of standard type}, $\stars$,
in the following way (which makes sense since $\srep$ is a bijection
into the space of all differential operators with smooth coefficients)
\BEQ {FlatStandardOrderedProd}
       \srep(F\stars G):= \srep(F)\srep(G)
\EEQ
and takes the following form:
\BEQ {standardflat}
       F\stars G = \sum_{r=0}^\infty
                   \frac{1}{r!}
                   \left(\frac{\hbar}{\im}\right)^r
     \frac{\partial^r F}{\partial p_{i_1} \cdots \partial p_{i_r}}
     \frac{\partial^r G}{\partial q^{i_1} \cdots \partial q^{i_r}} .
\EEQ
Due to (\ref{WeylundStandard}) we clearly have equivalence of $\stars$
and $\starw$
\BEQ {FlatWeylStdEquiv}
      F\stars G = N((N^{-1}F)\starw (N^{-1}G)) .
\EEQ
Writing $\pi$ for the canonical projection
$\mathbb R^{2n} \to \mathbb R^n: (q,p) \mapsto q$ and $i$ for the
canonical zero section $\mathbb R^n \to \mathbb R^{2n}: q \mapsto (q,0)$
we easily get the following useful formula for any smooth complex-valued
function $\psi$ on $\mathbb R^n$
\BEQ {isternpistern}
      \srep(F)\psi = i^*(F\stars (\pi^*\psi)) .
\EEQ
Considering $\hbar$ now as an additional variable on which the functions
depend polynomially we define the following differential operator
\BEQ {HomogenDerivDef}
     \mathcal H := p_i\frac{\partial}{\partial p_i} +
     \hbar\frac{\partial}{\partial \hbar} .
\EEQ
It is easy to see that both representations enjoy the following
{\em homogeneity property}:
\BEQ {FaltHomogenRep}
    \left[\hbar\frac{\partial}{\partial \hbar},
          \varrho_{\mbox{\rm \tiny S/W}}(F) \right]
    = \varrho_{\mbox{\rm \tiny S/W}}(\mathcal H F).
\EEQ
Physically this means that the operator corresponding to the
momentum component $p_l$ has also the physical dimension
of a momentum which is equal to the dimension of $\hbar$ divided
by length (the dimension of $q^l$). As a consequence, the two
star products are also homogeneous in the sense that the map
$\mathcal H$ is a derivation:
\BEQ {FlatHDeriv}
    \mathcal H (F*_{\mbox{\rm \tiny S/W}} G) =
                    (\mathcal H F)*_{\mbox{\rm \tiny S/W}} G
                  + F*_{\mbox{\rm \tiny S/W}} (\mathcal H G).
\EEQ

Both star products are associative and deform the pointwise
multiplication such that the component of the commutator which
is first order in $\hbar$ equals ${\im}$ times the canonical
Poisson bracket. Moreover, as can easily be seen by the
formulas (\ref{Weylflat}) and (\ref{standardflat}) both
star products are {\em bidifferential} in each order of
$\hbar$. They are even of {\em Vey type} which means
that in order $\hbar^r$ the corresponding bidifferential
operator is of order $r$ in each argument. One might be tempted to
think that any ordering prescription may give rise to a
reasonable star product, but the following example indicates
that one may lose the property that the star product
be bidifferential in each order of $\hbar$:

For $\mathbb R^2$ define the following modified ordering prescription
for a positive real number $s$:
\BEQ {pervers}
 \varrho_{\mbox{\rm\tiny perv}}((q,p)\mapsto f(q)p^k)
        := \left\{ \begin{array}{ccc}
                   \srep((q,p)\mapsto f(q)p^k)
                   & \mbox{ if } & k \neq 2 \\
                   \srep((q,p)\mapsto f(q)p^2) +
                   s\frac{\hbar}{\im}\srep((q,p)\mapsto f(q)p)
                   & \mbox{ if } & k=2
                   \end{array}
           \right.
\EEQ
A lengthy, but straight forward computation using formal symbols
$F(q,p):=e^{\alpha q+\beta p}$ where $\beta$ is considered as a formal
parameter shows that in the corresponding star product
(which is well-defined on all smooth functions polynomial in the momenta),
i.~e.~$\varrho_{\mbox{\rm\tiny perv}}(F)
        \varrho_{\mbox{\rm\tiny perv}}(G)
 = \varrho_{\mbox{\rm\tiny perv}}(F*_{\mbox{\rm\tiny perv}}G)$,
for each order in $\hbar$ there are infinitely many derivatives with
respect to $p$. Hence this kind of star product is in general {\em not}
extendable to the space of all formal power series in $\hbar$ with
coefficients in the smooth complex-valued functions on $\mathbb R^2$.

The aim is now to construct
and compute in more detail concrete star products on arbitrary cotangent
bundles $T^*Q$ and
possible representations as formal series of differential operators
on the formal series of smooth complex-valued functions on
$Q$. The important feature of these structures will be their
{\em homogeneity in the momenta} which is now defined using the
Liouville vector field $\xi$ on $T^*Q$ (whose flow consists in multiplying
the fibres by $e^t$) and which takes the familiar form
$p_i\frac{\partial}{\partial p_i}$ in a bundle chart.
It is interesting to note that the existence proof for star products
on arbitrary cotangent bundles $T^*Q$
by DeWilde and Lecomte \cite{DL83a} in 1983 is much easier than the general
proof thanks to this notion of homogeneity (which had earlier been used
by Cahen and Gutt for parallelizable manifolds, see \cite{CG82}):
By demanding the bidifferential operators
$M_r$ in the formal series of the star product of two smooth
complex-valued functions $f$ and $g$,
\BEQ {StarProduct}
    f\ast g = \sum_{r=0}^\infty \hbar^r M_r(f,g),
\EEQ
to be homogeneous of order $-r$ with respect to the Liouville field
DeWilde and Lecomte were able to show
that the usual obstructions in the third de Rham cohomology which a
priori occur when constructing the $M_r$ by induction
(see e.~g.~\cite{DL88}) simply vanish due to the homogeneity requirement.

The immediate generalization of the standard representation
$\srep$ to an arbitrary cotangent bundle
$T^*Q$ proceeds as follows: take the space of sections
$\Gamma(\bigvee^kTQ)$ of the $k$-fold symmetrized tangent bundle
and consider the canonical linear injection
$\widehat{~}:\Gamma(\bigvee^kTQ)\rightarrow
C^\infty(T^*Q):
T\mapsto \widehat{T}:\alpha\mapsto\frac{1}{k!}T(\alpha,\ldots,\alpha)$.
We call the complexification of its image $C^\infty_{pp,k}(T^*Q)$, and the
direct sum
$\bigoplus_{k=0}^\infty C^\infty_{pp,k}(T^*Q)$
is denoted by $C^\infty_{pp}(T^*Q)$ which is the obvious analogue of the
smooth complex-valued functions polynomial in the momenta. Clearly
$\Lie_\xi F = kF$ iff $F$ in $C^\infty_{pp,k}(T^*Q)$. Moreover
choose a torsion-free
covariant derivative $\nabla_0$ in the tangent bundle of $Q$
and replace partial by covariant derivatives in (\ref{StandOrdFlat})
in the following manner: for $T=T_1+\cdots+T_N$ where
$T_k\in\Gamma(\bigvee^kTQ)$ ($0\leq k\leq N$)  we define for a smooth
complex-valued function $\psi$ on $Q$
\BEQ {StandOrdCov}
    \srep(\widehat{T})\psi : q\mapsto \sum_{k=0}^N
    \frac{1}{k!} \left(\frac{\hbar}{\im}\right)^k
    T_k^{i_1\cdots i_k} (q)
    i_s (\partial_{q^{i_1}}) \cdots i_s (\partial_{q^{i_k}})
    D_0^{(k)} \psi (q)
\EEQ
where $D_0^{(k)}$ is the $k$-fold symmetrized covariant
derivative, $i_s$ means symmetric substitution and the functions
$T_k^{i_1 \cdots i_k}$ are regarded as the components of the
contravariant symmetric tensor field $T_k$ on $Q$.
This will clearly induce a star product on
the commutative subalgebra $C^\infty_{pp}(T^*Q)$ of $C^\infty(T^*Q)$.
However, since the Christoffel symbols of the connection $\nabla_0$
modify the differential operators of flat space representation
by terms of lower order it is a priori questionable in view of the above
counterexample whether this star product is bidifferential at each order
of $\hbar$. In fact it turns out that it is bidifferential
(which seems to have been shown by hard
analytic techniques in the past, cf.~e.~g.~\cite{Pfl95}) and in
this paper we want to deal with this question in a more algebraic
manner. Our main results are the following:
\begin{itemize}
\item After giving a short generalizing review of the Fedosov
      construction on an arbitrary manifold in Section
      \ref {FedosovSec} (where we need
      the covariant derivative term in more generality and a fibrewise
      formulation of equivalence transformation)
      we then build up and compute the general Fedosov machinery for a
      star product of standard ordered $\stars$ and Weyl type $\starf$ on
      any cotangent
      bundle $T^*Q$ (Section \ref {WeylSec} and \ref {StandardSec}) 
      based on a rather natural, seemingly well-known lift
      of any torsion-free covariant derivative on $Q$ to a symplectic
      torsion-free covariant derivative on $T^*Q$ which is homogeneous
      with respect to the Liouville vector field on $T^*Q$ (see Appendix
      \ref {ConnApp} for a description). It turns out
      that the corrections to the covariant derivative needed to
      construct the Fedosov derivative are ``classical'' in the sense
      that they do {\em not} depend on the formal parameter. The standard
      ordered type fibrewise star product is constructed by using the
      duality between the horizontal and the vertical subbundle of the
      tangent bundle of $T^*Q$. Then we can prove by constructing an
      equivalence transformation in Section \ref {EquiSec} that the 
      Fedosov star product of standard ordered type $\stars$ is equivalent to
      the Fedosov star product of Weyl type $\starf$.

\item According to the general Fedosov philosophy ``Whatever you plan
      to do on a symplectic manifold $M$, do it first fibrewise on the
      tangent spaces, and pull it then down to $M$ by means of a nice
      compatible Fedosov derivative'' we are then constructing a
      fibrewise standard representation analogous to the flat space
      formula (\ref{isternpistern}) in Section \ref {RepSec}. 
      As a surprise it turned out that our
      naively constructed Fedosov derivative $\mathcal D'$
      (\ref{DummyDDef}) (as a conjugate of the
      original Weyl type Fedosov derivative $\Dfed$ by means of
      the fibrewise
      analogue to the operator $N$, see (\ref{SFlat})) was {\em not}
      compatible with the fibrewise standard representation. Luckily,
      it could be modified by a fibrewise internal automorphism
      (see Thm.~\ref{hAutoTheo}) to
      render it compatible. As a result we obtain exactly the above standard
      representation (\ref{StandOrdCov}). Since it is easy to see that
      all `reasonable' star products constructed by a Fedosov type
      procedure automatically are bidifferential the construction shows
      a posteriori that the standard representation does not fall under
      the above-mentioned ``beasty'' class of ordering prescriptions.

\item Finally we derive a surprisingly simple analogue of the operator
      $N$ (cf.~(\ref{SFlat})) for any $T^*Q$ in Section \ref {WeylRepSec}: 
      it takes the form $N = \exp(\frac{\hbar}{2{\im}}\Delta)$ where 
      the second-order differential operator $\Delta$ takes the following 
      form in a bundle chart $(q,p)$:
      \BEQ {DerNeumeier}
           (\Delta F)(q,p) =
           \frac{\partial^2 F}{\partial q^i\partial p_i}(q,p)
           + \Gamma^i_{ik}(q)\frac{\partial F}{\partial p_k}(q,p)
           + p_r\Gamma^r_{ij}(q)
             \frac{\partial^2 F}{\partial p_i\partial p_j}(q,p)
           + \alpha_r(q)\frac{\partial F}{\partial p_r}(q,p)
      \EEQ
      where the $\Gamma^i_{jk}$ are the Christoffel symbols of the
      connection $\nabla_0$ and $\alpha$ is a particular choice of a 
      one-form on $Q$ such that $-d\alpha$ equals the trace of the 
      curvature tensor (see the Appendix for a Theorem). 
      In case $\nabla_0$ leaves invariant a volume on $Q$ (assumed to be 
      orientable) then $\alpha$ can be chosen to be zero. $N$ can now be
      used as an equivalence transformation from the star product of
      standard ordered type $\stars$ to another star product $\starw$
      (\ref{TrueWeylProdDef})
      which is of Weyl type but which turns out to be
      {\em different} from the Fedosov star product of Weyl type,
      $\starf$!

      As an example, we rederive the star product on the cotangent bundle
      of an arbitrary Lie group constructed by means of the standard
      torsion-free left-invariant `half commutator' connection first
      given by Gutt (see \cite{Gut83}) and give an explicit closed
      formula of a star product of standard ordered type in Section 
      \ref {LieSec}.
\end{itemize}

\vspace{0.5cm}

\noindent {\bf Convention:} In what follows $\hbar$ will always denote
a real number whereas $\lambda$ will denote a formal parameter
which in converging situations may be substituted by $\hbar$ and
is considered to be real, i.~e.~$\cc \lambda := \lambda$.

\section {Fedosov Derivations and Fedosov-Taylor series}
\label {FedosovSec}

In this rather technical section we shall revisite Fedosov's
construction of star products in a slightly more general
context. The notation is mainly the same as in
Fedosov's book \cite {Fed96} and in \cite {BW96a}.
Let $M$ be a smooth manifold and define
\BEQ {WeylAlgDef}
    \WL (M) := \left({\mathsf X}_{s=0}^\infty \mathbb C\left(
               \Gamma\left(\mbox{$\bigvee$}^s T^*M
               \otimes \mbox{$\bigwedge$}
               T^*M\right)\right)\right)[[\lambda]] .
\EEQ
If there is no possibility for confusion we simply write
$\WL$ and denote by $\WL^k$ the elements of antisymmetric degree $k$ and
set $\W := \WL^0$. For two elements $a, b \in \WL$ we define their
pointwise product denoted by $\mu (a \otimes b) = ab$
by the symmetric $\vee$-product in the first factor
and the antisymmetric $\wedge$-product in the second factor. Then
the degree-maps $\degs$ and $\dega$ with respect to the
symmetric and antisymmetric degree are derivations of this product.
Therefore we shall call $\WL$ a formally
$\mathbb Z \times \mathbb Z$-graded
algebra with respect to the symmetric and antisymmetric degree.
Moreover $(\WL, \mu)$ is supercommutative with respect to the
antisymmetric degree.
For a vector field $X$ we define the symmetric substitution (insertion) 
$i_s(X)$ and the antisymmetric substitution $i_a (X)$ which are
superderivations of symmetric degree $-1$ resp. $0$ and
antisymmetric degree $0$ resp. $-1$. Following Fedosov we define
\BEQ {deltaDef}
    \delta := (1 \otimes dx^i) i_s (\partial_{x^i})
    \quad \mbox { and } \quad
    \delta^* := (dx^i \otimes 1) i_a (\partial_{x^i})
\EEQ
where $x^1, \ldots, x^n$ are local coordinates for $M$
and for $a \in \WL$ with $\degs a = ka$ and $\dega a = la$
we define
\BEQ {deltaInvDef}
    \delta^{-1}a := \left\{
    \begin {array} {cl}
    \frac{1}{k+l} \delta^*a & \mbox { if } k+l \ne 0 \\
    0 & \mbox { if } k+l=0
    \end {array}
    \right.
\EEQ
and extend $\delta^{-1}$ by linearity. Clearly
$\delta^2 = {\delta^*}^2 = 0$. Moreover we denote by
$\sigma: \WL \to C^\infty (M)[[\lambda]]$ the projection onto the
part of symmetric and antisymmetric degree $0$.
Then one has the following `Hodge-decomposition' for any
$a \in \WL$ (see e.~g. \cite [eq. 2.8.] {Fed94}):
\BEQ {HodgeDecomp}
    a = \delta \delta^{-1} a + \delta^{-1} \delta a + \sigma (a)
\EEQ
Now we consider a fibrewise associative deformation $\circ$
of the pointwise product which should have the form
\BEQ {FibProd}
    a \circ b = ab + \sum_{r=1}^\infty \lambda^r
                \mathcal M_r (a, b)
\EEQ
where
$\mathcal M_r (a, b) = M^{i_1\ldots i_r j_1 \ldots j_r}_r
i_s (\partial_{x^{i_1}}) \cdots i_s (\partial_{x^{i_r}}) a
i_s (\partial_{x^{j_1}}) \cdots i_s (\partial_{x^{j_r}}) b$ and the
$M_r^{i_1 \ldots i_r j_1 \ldots j_r}$ are the coefficients of a
tensor field totally symmetric in $i_1, \ldots, i_r$ and
$j_1, \ldots, j_r$ separately.
Moreover we define $\dega$-graded supercommutators with respect to
$\circ$ and set $\ad (a) b := [a, b]$. If not all $\mathcal M_r = 0$
for $r \ge 1$ then $\degs$ is no longer a derivation of the deformed
product $\circ$ but $\Deg := \degs + 2 \degh$
is still a derivation and hence the algebra
$(\WL, \circ)$ is formally $\Deg$-graded where
$\degh := \lambda \partial_\lambda$.
We shall refer to this degree as total degree.
We shall treat the non-deformed case separately at the end of this
section and first remember the following two theorems which can
be proved completely analogously to Fedosov's original theorems in
\cite[Theorem 3.2, 3.3]{Fed94}:
\begin {THEOREM} \label {GenFDerivTheo}
Let $T^{(0)}: \WL \to \WL$ be a superderivation of $\circ$ of
antisymmetric degree $1$ and total degree $0$ such that
$[\delta, T^{(0)}] = 0$, and
$\left(T^{(0)}\right)^2 = \frac{1}{2} [T^{(0)}, T^{(0)}] =
\frac{\im}{\lambda} \ad (T)$ with some $T \in \WL^2$ of total degree $2$,
and let $T$ satisfy $\delta T = 0 = T^{(0)} T$.
Then there exists a unique element $r \in \WL^1$ such that
\BEQ {Genr}
    \delta r = T + T^{(0)} r + \frac{\im}{\lambda} r \circ r
    \quad \mbox { and } \quad
    \delta^{-1} r = 0 .
\EEQ
Moreover $r = \sum_{k=3}^\infty r^{(k)}$ with
$\Deg r^{(k)} = k r^{(k)}$ satisfies the recursion formulas
\BEQ {GenrRecus}
    \begin {array} {c}
    r^{(3)} = \delta^{-1} T \\
    \displaystyle
    r^{(k+3)} = \delta^{-1} \left( T^{(0)} r^{(k+2)} +
                \frac{\im}{\lambda}
                \sum_{l=1}^{k-1} r^{(l+2)} \circ r^{(k-l+2)} \right) .
    \end {array}
\EEQ
In this case the Fedosov derivation
\BEQ {GenFedDerivDef}
    \mathcal D := - \delta + T^{(0)} + \frac{\im}{\lambda} \ad (r)
\EEQ
is a superderivation of antisymmetric degree $1$ and
has square zero: $\mathcal D^2 = 0$.
\end {THEOREM}
\begin {THEOREM} \label {GenFTaylorTheo}
Let $L = -\delta + T: \WL \to \WL$ be a $\mathbb
C[[\lambda]]$-linear map of antisymmetric degree $1$ with square
zero $L^2 = 0$ such that $T$ does not decrease the total degree.
\begin {enumerate}
\item Then for any $f \in C^\infty (M)[[\lambda]]$ there exists a
      unique element $\tau_L (f) \in \ker L \cap \W$ such that
      \BEQ {sigmaTaylorL}
          \sigma (\tau_L (f)) = f
      \EEQ
      and $\tau_L : C^\infty (M)[[\lambda]] \to \W$
      is $\mathbb C[[\lambda]]$-linear and refered to as the
      Fedosov-Taylor series corresponding to $L$.

\item If in addition $T = \sum_{k=0}^\infty T^{(k)}$ such that
      $T^{(k)}$ is homogeneous of total degree $k$ then
      for $f \in C^\infty (M)$ we have
      $\tau_L (f) = \sum_{k=0}^\infty \tau_L (f)^{(k)}$ where
      $\Deg\tau_L (f)^{(k)} = k \tau_L (f)^{(k)}$ which can
      be obtained by the following recursion formula
      \BEQ {GenTaylorRecurs}
          \begin {array} {c}
              \tau_L (f)^{(0)} = f \\
              \displaystyle
              \tau_L (f)^{(k+1)} =
              \delta^{-1} \sum_{l=0}^k T^{(l)} \tau_L (f)^{(k-l)}.
          \end {array}
      \EEQ

\item If $L = \mathcal D$ is a $\circ$-superderivation of
      antisymmetric degree $1$ as constructed in theorem
      \ref {GenFDerivTheo} then $\ker \mathcal D \cap \W$ is 
      a $\circ$-subalgebra and a new (eventually deformed) 
      associative product $*_{\mathcal D}$ for $C^\infty (M)[[\lambda]]$ 
      is defined by pull-back of $\circ$ via $\tau_{\mathcal D}$.
\end {enumerate}
\end {THEOREM}
Let $\circ'$ be another fibrewise product for $\WL$ and
$\mathcal S = \id + \sum_{r=1}^\infty \lambda^r \mathcal S_r$
where $\mathcal S_r: \WL \to \WL$ is a map of the form
$\mathcal S_r = \sum_{l=1}^{2r} \frac{1}{l!} S^{i_1 \cdots i_l}
i_s (\partial_{x^{i_1}}) \cdots i_s (\partial_{x^{i_l}})$ where
$S^{i_1 \cdots i_l}$ are the compoments of a symmetric tensor
field. If in addition
\[
    \mathcal S (a \circ b) = (\mathcal Sa) \circ' (\mathcal Sb)
\]
for all $a, b \in \WL$ then $\mathcal S$ is called a fibrewise equivalence
transformation between $\circ$ and $\circ'$. Note that $\mathcal S$ is
clearly invertible and $\mathcal S^{-1}$ is a fibrewise equivalence
transformation between $\circ'$ and $\circ$.
If $\mathcal D$ is a $\circ$-superderivation as constructed in
theorem \ref {GenFDerivTheo} and $\tau_{\mathcal D}$ its corresponding
Fedosov-Taylor series and $*_{\mathcal D}$ the induced associative
product for $C^\infty (M)[[\lambda]]$ then
$\mathcal D' := \mathcal S \mathcal D \mathcal S^{-1}$ is a
$\circ'$-superderivation of antisymmetric degree $1$ of the form
$\mathcal D' = -\delta + T'$
satisfying the conditions of part one of the preceding theorem.
Hence there exists a corresponding Fedosov-Taylor series
$\tau_{\mathcal D'}$ which induces an associative product
$*_{\mathcal D'}$ on $C^\infty (M)[[\lambda]]$. Then $*_{\mathcal D}$
and $*_{\mathcal D'}$ turn out to be equivalent too:
\begin {PROPOSITION} \label {EquivalenceTrafoProp}
With the notation form above we define the map
$T: C^\infty (M)[[\lambda]] \to C^\infty (M)[[\lambda]]$
\BEQ {EquivalenceTrafoDef}
    Tf := \sigma (\mathcal S \tau_{\mathcal D} (f))
\EEQ
for $f \in C^\infty (M)[[\lambda]]$ which is an
$\mathbb C[[\lambda]]$-linear equivalence transformation
between $*_{\mathcal D}$ and $*_{\mathcal D'}$, 
i.~e. $T ( f *_{\mathcal D} g) = (Tf) *_{\mathcal D'} (Tg)$ for all
$f, g \in C^\infty (M)[[\lambda]]$ with
inverse $T^{-1} f = \sigma (\mathcal S^{-1} \tau_{\mathcal D'} (f))$.
\end {PROPOSITION}
\begin {PROOF}
This is a straight forward computation observing
$\mathcal D' \mathcal S \tau_{\mathcal D} (f) = 0$ and applying the
last theorem.
\end {PROOF}

At last we shall discuss the classical case with the undeformed
product $\mu$. First we restrict our considerations to the
classical part $\WL_{\rm cl}$ of $\WL$ which are just those
elements without any positive $\lambda$-powers. Next we
consider a torsion-free connection $\nabla$ for $M$ and define the
map (using the same symbol as for the connection)
\BEQ {NablaDef}
    \nabla := (1 \otimes dx^i) \nabla_{\partial_{x^i}}
\EEQ
where $\nabla_{\partial_{x^i}}$ denotes the covariant derivative
with respect to $\partial_{x^i}$. Then clearly $\nabla$ is
globally defined and a superderivation of $\mu$ which leaves
$\WL_{\rm cl}$ invariant. Moreover we consider
\BEQ {WLXDef}
    \WLX (M) := \left({\mathsf X}_{s=0}^\infty \mathbb C\left(
                \Gamma\left(\mbox{$\bigvee$}^s T^*M
                \otimes \mbox{$\bigwedge$}
                T^*M \otimes TM\right)\right)\right)[[\lambda]]
\EEQ
and define $\WLX_{\rm cl}$ analogously. For $\varrho \in \WLX$
we define the symmetric substitution $i_s (\varrho)$ by inserting the
vector part of $\varrho$ symmetrically and multiplying the form
part of $\varrho$ by $\mu$ from the left.
Then we have $\nabla^2 = -i_s (R)$ where $R$ is the curvature
tensor viewed as element of antisymmetric degree $2$ in $\WLX$.
The classical analogue to the Fedosov derivation is decribed by the
following theorem which is due to Emmrich and
Weinstein \cite [Theorem 1] {EW93}:
\begin {THEOREM} \label {ClassFedDerivTheo}
Let $\nabla$ be defined as in (\ref {NablaDef}) then $\nabla$ is a
superderivation of the undeformed product $\mu$ and there exists a
uniquely determined element $\varrho_0 \in \WLX_{\rm cl}$ of
antisymmetric degree $1$ such that $\delta^{-1} \varrho_0 = 0$
and such that the classical Fedosov derivation
\BEQ {ClassFDerivDef}
    \mathcal D_0 = -\delta + \nabla + i_s (\varrho_0)
\EEQ
has square zero: $\mathcal D_0^2 = 0$. Note that
$\mathcal D_0 \varrho_0 = R_0$ and $\cc{\varrho_0} = \varrho_0$ as well as
$\cc {\mathcal D_0 a} = \mathcal D_0 \cc a$ for all $a \in \WL$.
\end {THEOREM}
Furthermore in \cite [Theorem 3 and 6] {EW93} was shown
using analytical techniques
that in this case the corresponding Fedosov-Taylor series $\tau_0$
is just the formal Taylor series with respect to the connection.
We shall give here another more algebraic proof of this result:
\begin {THEOREM} \label {ClassFedTaylorTheo}
Let $\tau_0$ be the Fedosov-Taylor series of $\mathcal D_0$
according to theorem \ref {GenFTaylorTheo}. Then for
$f \in C^\infty (M)$ we have
\BEQ {ClassTau}
    \tau_0 (f) = e^{D} f
\EEQ
where $D = dx^i \vee \nabla_{\partial_{x^i}}$ and hence
$\tau_0 (f)$ is the formal Taylor series with repect to the
connection $\nabla$.
\end {THEOREM}
\begin {PROOF}
Since $\mathcal D_0$ satisfies all conditions for theorem
\ref {GenFTaylorTheo} we only have to show that $e^D f$ satisfies
the recursion formula (\ref {GenTaylorRecurs}). Note that in this
particular case the total degree and the symmetric degree
coincide. First we observe $D = [\delta^*, \nabla]$ and hence
applying $\delta$ and $\delta^*$ to (\ref {GenTaylorRecurs}) implies
$0 = - (\delta^*\delta + \delta \delta^*) \tau_0 (f)^{(k+1)}
     + D\tau_0 (f)^{(k)} + \sum_{l=0}^k \delta^*
     \left(i_s (\varrho_0^{(l+1)}) \tau_0 (f)^{(k-l)}\right)$.
Now each term in the last sum vanishes identically due to
$\delta^* \varrho_0 = 0$ and $\dega \tau_0 (f) = 0$ hence
$(k+1) \tau_0(f)^{(k+1)} = D\tau_0 (f)^{(k)}$ since
$\delta^*\delta + \delta \delta^* = \degs + \dega$. Then
(\ref{ClassTau}) follows directly by induction on the symmetric
degree $k$.
\end {PROOF}

\section {Homogeneous Fedosov star product of Weyl type}
\label {WeylSec}

Let $\pi: T^*Q \to Q$ be the cotangent bundle of a differentiable,
$n$-dimensional manifold $Q$ and let $\theta_0$ be the canonical
one-form, $\omega_0 := -d\theta_0$ the canonical symplectic form and
$\xi$ defined by $i_\xi \omega_0 = -\theta_0$ the canonical (Liouville)
vector field on $T^*Q$. Moreover let $i: Q \to T^*Q$ be the embedding of
$Q$ in $T^*Q$ as zero section. We consider now for $\WL$ of $T^*Q$
the fibrewise Weyl product defined for $a, b \in \WL$ by
\BEQ {FibWeylProdDef}
    a \fed b := \mu \circ e^{\frac{\im\lambda}{2} \Lambda^{kl}
                 i_s (\partial_{x^k}) \otimes i_s (\partial_{x^l})}
                 a \otimes b
\EEQ
where $\mu (a\otimes b) = ab$ is the fibrewise product in $\WL$
and $\Lambda^{kl}$ are the components of
the canonical Poisson tensor with respect to some coordinates
$x^1, \ldots, x^{2n}$ of $T^*Q$. It will be advantageous for
calculations to use local bundle (Darboux) coordinates $q^1, \ldots, q^n$,
$p_1, \ldots, p_n$ induced by coordinates $q^1, \ldots, q^n$ on $Q$
such that $p_1, \ldots, p_n$ are the conjugate momenta to the
$q^1, \ldots, q^n$. In the following
$q^1, \ldots, q^n, p_1, \ldots p_n$  shall always denote such a
local bundle (Darboux) chart. Obviously $\fed$ is an associative
deformation of $\mu$ of the form as in (\ref {FibProd}) and hence we
can apply all results of section \ref {FedosovSec} to this
particular situation. We denote by $\Lie_\xi$ the Lie derivative
with respect to the canonical vector field $\xi$ and define the
`homogeneity derivation'
\BEQ {HDef}
    \mathcal H := \Lie_\xi + \degh = \Lie_\xi +
    \lambda \frac{\partial}{\partial\lambda} .
\EEQ
Note that $\mathcal H$ is $\mathbb C$-linear but not
$\mathbb C[[\lambda]]$-linear. Then the fibrewise Weyl product is
homogeneous which can be proved as in the case of $\mathbb R^{2n}$:
\begin {LEMMA}
Let $a, b \in \WL$
then $\mathcal H$ is a (super-)derivation of $\fed$ of
antisymmetric and total degree $0$:
\BEQ {HDeriv}
    \mathcal H (a \fed b) = \mathcal H a \fed b + a \fed \mathcal H b
\EEQ
Moreover we have
$\left[ \Lie_\xi, \delta \right] =
\left[ \Lie_\xi, \delta^* \right] =
\left[ \Lie_\xi, \delta^{-1} \right] = 0$ and
$\left[ \mathcal H, \delta \right] =
\left[ \mathcal H, \delta^* \right] =
\left[ \mathcal H, \delta^{-1} \right] = 0$.
\end {LEMMA}

According to Fedosov's construction of a star product one needs
a torsion-free and symplectic connection $\nabla$ for $T^*Q$ and the
map $\nabla: \WL \to \WL$ defined as in (\ref {NablaDef}).
If the connection is symplectic $\nabla$ turns out to be a
superderivation of antisymmetric degree $1$ and symmetric
and total degree $0$ of the fibrewise Weyl product $\fed$. Moreover
$[\delta, \nabla] = 0$ and $2\nabla^2 = [\nabla, \nabla]$ turns
out to be an inner superderivation
\BEQ {NablaSquare}
    \nabla^2 = \frac{\im}{\lambda} \adfed (R)
\EEQ
where $R := \frac{1}{4} \omega_{it} R^t_{jkl}
dx^i \vee dx^j \otimes dx^k \wedge dx^l \in \WL^2$ involves the
curvature of the connection.
Moreover one has $\delta R = 0 = \nabla R$ as a consequence of
the Bianchi identities. Hence the map $\nabla$ as
term of total degree $0$ satisfies all conditions of theorem
\ref {GenFDerivTheo} and hence there is a unique element
$\rfed \in \WL^1$ such that $\delta \rfed =
R + \nabla\rfed + \frac{\im}{\lambda} \rfed \fed \rfed$
and $\delta^{-1} \rfed = 0$. Moreover the Fedosov derivation
\BEQ {FDerivWeylDef}
    \Dfed := - \delta + \nabla
              + \frac{\im}{\lambda} \adfed (\rfed)
\EEQ
has square $0$. Let $\taufed$ be the corresponding
Fedosov-Taylor series. Then Fedosov has shown that
\BEQ {WeylProdDef}
    f \starf g := \sigma (\taufed (f) \fed \taufed (g))
\EEQ
defines a star product \cite [eq. 3.14]{Fed94} which is of Weyl
type (see e.~g.~\cite [Lemma 3.3] {BW96a}).

In the particular case of a cotangent
bundle we consider homogeneous, symplectic and torsion-free
connections (see definition \ref {HomConnDef}):
\begin {LEMMA}
If the symplectic and torsion-free connection $\nabla$ on $T^*Q$ is
in addition homogeneous then
\BEQ {LieXiNabla}
    [\Lie_\xi , \nabla] = [\mathcal H, \nabla] = 0
\EEQ
\BEQ {LieXiR}
    \mathcal HR = \Lie_\xi R = R .
\EEQ
\end {LEMMA}
\begin {THEOREM} \label {WeylProTheo}
Let $\nabla$ be defined as above with the additional property that
the connection is homogeneous. Then
\BEQ {LieXirWeyl}
    \Lie_\xi \rfed = \rfed = \mathcal H \rfed
    \quad \mbox { and }
    \quad
    \frac{\partial}{\partial \lambda} \rfed = 0 .
\EEQ
Moreover the Fedosov derivation $\Dfed$ (super-)commutes with
$\mathcal H$
\BEQ {HomFDeriv}
    \left[\mathcal H, \Dfed \right] = 0
\EEQ
and $\rfed$ satisfies the following simpler recursion formulas
\BEQ {rSimpleRec}
    \begin {array} {c}
    \rfed^{(3)} = \delta^{-1} R \\
    \displaystyle
    \rfed^{(k+3)} =
    \delta^{-1} \left(\nabla \rfed^{(k+2)} -
                \frac{1}{2} \sum_{l=1}^{k-1} \left\{
                \rfed^{(l+2)}, \rfed^{(k-l+2)}
                \right\}_{\mbox {\rm fib}} \right)
    \end {array}
\EEQ
where $\{\cdot, \cdot\}_{\rm fib}$ denotes the fibrewise
Poisson bracket in $\WL$. Moreover the corresponding Fedosov-Taylor
series $\taufed$ commutes with $\mathcal H$
\BEQ {HomTauWeyl}
    \mathcal H \taufed (f) = \taufed (\mathcal H f)
\EEQ
and satisfies the usual recursion formulas for $f \in C^\infty (T^*Q)$ 
with respect to the total degree 
\BEQ {FedTaylorRecurs}
    \begin {array} {c}
        \taufed (f)^{(0)} = f \\
        \displaystyle
        \taufed (f)^{(k+1)} = \delta^{-1} \left(
        \nabla \taufed (f)^{(k)} + 
        \frac{\im}{\lambda} \sum_{l=1}^{k-1} 
        \adfed \left(\rfed^{(l+2)} \right) \taufed (f)^{(k-l)} \right)
    \end {array}
\EEQ            
analogously to (\ref {GenTaylorRecurs}). The Fedosov star product 
(\ref {WeylProdDef}) is homogeneous, i.~e. for 
$f, g \in C^\infty (T^*Q)[[\lambda]]$ 
\BEQ {HomWeylStarProd}
    \mathcal H (f \starf g)
    = (\mathcal H f) \starf g + f \starf (\mathcal H g) .
\EEQ
\end {THEOREM}
\begin {PROOF}
The fact $\mathcal H \rfed = \rfed$ is proved by
induction using the recursion formulas for $\rfed$. Then
(\ref {HomFDeriv}), (\ref {HomTauWeyl}) and
(\ref {HomWeylStarProd}) easily follow and (\ref {FedTaylorRecurs}) 
follows directly from (\ref {GenTaylorRecurs}). Since
$\mathcal H \rfed = \rfed$ the section $\rfed$ can depend at
most linearly on $\lambda$ but since it has to depend on even
powers of $\lambda$ only (see \cite [Lemma 3.3]{BW96a}) it has to be
independent of $\lambda$ at all. Then the recursion
formulas (\ref {rSimpleRec}) follow by induction.
\end {PROOF}

We shall refer to $\starf$ as the {\em homogeneous Fedosov star product
of Weyl type} induced by the homogeneous connection $\nabla$.
Using this theorem we find several corollaries. The first one is
originally due to DeWilde and Lecomte \cite [Proposition 4.1]{DL83a}:
\begin {COROLLARY}
On every cotangent bundle $T^*Q$ there exists a homogeneous star
product of Weyl type.
\end {COROLLARY}
\begin {COROLLARY}
Let $f \in C^\infty_{pp, k} (T^*Q)$ then $\taufed (f)$ contains
only even powers of $\lambda$ up to order $k$.
\end {COROLLARY}
\begin {PROOF}
This easily follows from (\ref {HomTauWeyl}) and \cite [Lemma 3.3]{BW96a}.
\end {PROOF}
\begin {COROLLARY}
The Fedosov-Taylor series $\taufed$ satisfies
\BEQ {taupiCom}
    \taufed \circ \pi^* = \pi^* \circ i^* \circ \taufed \circ \pi^* .
\EEQ
\end {COROLLARY}

General properties of homogeneous star products are decribed in the
following proposition:
\begin {PROPOSITION} \label {AnalProp}
Let $*$ be a homogeneous star product for $T^*Q$.
\begin {enumerate}
\item The functions polynomial in $\lambda$ and in the momenta
      $C^\infty_{pp} (T^*Q)[\lambda]$ are
      a $\mathbb C[\lambda]$-submodule of
      $C^\infty (T^*Q)[[\lambda]]$ with respect to $*$
      and hence $f*g$ trivially converges for all
      $\lambda = \hbar \in \mathbb R$ and
      $f, g \in C^\infty_{pp} (T^*Q)[\lambda]$.
\item Let $U \subseteq Q$ be a domain of a chart. Then any
      $f \in C^\infty_{pp} (T^*U)[\lambda]$ can be written as finite sum
      of star products of functions in $C^\infty_{pp,0} (T^*U)[\lambda]$
      and $C^\infty_{pp,1} (T^*U)[\lambda]$.
\item The vector space $C^\omega_{p} (T^*Q)[[\lambda]]$ of
      formal power series with coefficients in the functions which
      are analytic in the fibre variables are a
      $\mathbb C[[\lambda]]$-submodule of
      $C^\infty (T^*Q)[[\lambda]]$.
\end {enumerate}
\end {PROPOSITION}
\begin {PROOF}
The first part is obvious using lemma \ref {WichtigLem}
and for the third part one observes that the homogeneity implies that
the coefficient functions of the bidifferential
operators in the star product are polynomial in the momenta.
For the second part note that $f \in C^\infty_{pp, k+1} (T^*U)$ takes the
form $f (q, p) = \frac{1}{(k+1)!} F^{i_1 \cdots i_{k+1}}(q)
p_{i_1} \cdots p_{i_{k+1}}$. Subtracting
$\frac{1}{(k+1)!} (\pi^* F^{i_1 \cdots i_{k+1}})
* p_{i_1} * \cdots * p_{i_{k+1}}$ results in a finite sum of functions of
degree smaller than $k+1$ with respect to $\Lie_\xi$ and polynomial in
$\lambda$ proving the obvious induction on $k$.
\end {PROOF}

Now we consider the Fedosov derivation $\Dfed$ and the section $\rfed$
more closely. First we notice that $\delta$, $\delta^*$ and
$\delta^{-1}$ satisfy the following relations
\BEQ {deltaPi}
    \delta \pi^* = \pi^* \delta_0
    \quad
    \delta^* \pi^* = \pi^* \delta^*_0
    \quad
    \delta^{-1} \pi^* = \pi^* \delta^{-1}_0
\EEQ
where $\delta_0$, $\delta^*_0$ and $\delta^{-1}_0$ are the
corresponding maps on $Q$ defined analogously to $\delta$,
$\delta^*$ and $\delta^{-1}$. Moreover for a homogeneous connection
we get
\BEQ {NablaPi}
    \nabla \pi^* = \pi^* \nabla_0
\EEQ
where $\nabla_0$ is the corresponding map on $Q$ defined by the
induced connection $\nabla_0$ on $Q$
(see definition \ref {InducedConnDef}). By direct calculation we get
the following proposition:
\begin {PROPOSITION}
Let $\Dfed$ and $\rfed$ be given as in (\ref {FDerivWeylDef}).
Then there exists a unique element
$\varrho \in \WLX_{\rm cl}(Q)$ of antisymmetric degree $1$ such that
\BEQ {DweylPi}
    \Dfed \pi^* = \pi^* \mathcal D
\EEQ
where $\mathcal D = -\delta_0 + \nabla_0 + i_s (\varrho)$ and clearly
$\mathcal D^2 =0$. In local coordinates the element $\varrho$ 
takes the following form
\BEQ {ClassrhoLocal}
    \varrho = i^*(i_s (\partial_{p_i})\rfed) \otimes \partial_{q^i}
\EEQ
and we have $\delta^*_0 \varrho = 0$ iff $i_a (X) \rfed = 0$ for all
vertical vector fields $X \in \Gamma (T(T^*Q))$.
\end {PROPOSITION}
With other words, if $i_a (X) \rfed = 0$ for all vertical vector
fields then $\mathcal D$ would coincide with the
map $\mathcal D_0$ and $\varrho$ would coincide with $\varrho_0$
as in theorem \ref {ClassFedDerivTheo} applied for $M = Q$.
In the following we shall prove that for any torsion-free connection
on $Q$ there is indeed a canonical choice for a homogeneous connection
on $T^*Q$ such that this is the case:
\begin {PROPOSITION}
Consider a torsion-free connection on $Q$ and the map $\nabla_0$
defined as in (\ref {NablaDef}) and let $\nabla$ be a homogeneous,
symplectic and torsion-free connection for $T^*Q$ such that
$\nabla \pi^* = \pi^* \nabla_0$ and let $\rfed$
be the corresponding element in $\WL^1$. Then
$i_a (X) \rfed^{(3)} = 0$
for every vertical vector field $X \in \Gamma (T(T^*Q))$ where
$\rfed^{(3)}$ is the term of total degree $3$ in $\rfed$
iff the connection $\nabla$ coincides with the lifted connection
$\nabla^0$ defined as in (\ref {LiftedConnDef}). Moreover in this case
\BEQ {iaXrWeyl}
    i_a (X) \rfed = 0.
\EEQ
\end {PROPOSITION}
\begin {PROOF}
The fact $i_a (X) \rfed^{(3)} = 0$ for $X$ vertical iff
$\nabla = \nabla^0$ follows from proposition
\ref {iaXDeltaInvRProp}. Then (\ref {iaXrWeyl}) follows by a lengthy
but straight forward induction on the total degree using the
recursion formulas (\ref {rSimpleRec}).
\end {PROOF}

\begin {COROLLARY} \label {TauWeylPiCor}
Let $\nabla_0$ be a torsion-free connection for $Q$ and $\nabla^0$
the corresponding homogeneous, symplectic, and torsion-free
connection for $T^*Q$. Then the corresponding Fedosov derivation
$\Dfed$ satisfies
\BEQ {DweylPiDNull}
    \Dfed \pi^* = \pi^* \mathcal D_0
\EEQ
and hence the Fedosov-Taylor series of the pull-back of functions
$\chi \in C^\infty (Q)[[\lambda]]$ is just the pull-back of
their Taylor series with respect to $\nabla_0$
\BEQ {TauPi}
    \taufed (\pi^*\chi) = \pi^* \tau_0 (\chi) .
\EEQ
\end {COROLLARY}

At last in this section we shall discuss the classical limit of
the Fedosov derivation $\Dfed$ and the Fedosov-Taylor series
$\taufed$. We define $\Dfedcl$ by setting $\lambda = 0$ in
$\Dfed$ which is well-defined since in any case $\adfed (\rfed)$
starts with $\im\lambda \{\rfed, \cdot \}_{\rm fib}$ and analogously
we define $\taufedcl$. Then clearly for $f, g \in C^\infty (T^*Q)$
\BEQ {DweylClassic}
    (\Dfedcl)^2 = 0
    \quad
    \mbox { and }
    \quad
    \Dfedcl \taufedcl (f) = 0
    \quad
    \mbox { and }
    \quad
    \taufedcl (\{f, g\}) =
    \{\taufedcl (f), \taufedcl (g)\}_{\rm fib}
\EEQ
which is also true for arbitrary symplectic manifolds
(compare \cite [Sec. 8] {EW93}). Now we
concentrate again on the particular case where we have chosen a
torsion-free connection $\nabla_0$ for $Q$ and the corresponding
homogeneous, symplectic, and torsion-free connection $\nabla^0$
for $T^*Q$. Then for any vertical vector field
$X \in \Gamma (T(T^*Q))$ we have
\BEQ {iaXDweyl}
    i_a (X) \Dfed + \Dfed i_a (X) =
    \left[ i_a (X), \Dfed\right] =
    \Lie_X - i_s (X)
    - (dx^i \otimes 1) i_s (\nabla^0_{\partial_{x^i}} X)
\EEQ
which is proved by direct calculation using $i_a (X) \rfed = 0$.
Note that this equation is also true if $\Dfed$ is replaced by
$\Dfedcl$ since the right hand side is obviously independent of
$\lambda$.
\begin {THEOREM} \label {ClassFTTheo}
Let $\nabla_0$ be a torsion-free connection for $Q$ and
$\nabla^0$ the corresponding homogeneous, symplectic and torsion-free
connection for $T^*Q$. Moreover let $\taufed$ be the corresponding
Fedosov-Taylor series and $\taufedcl$ the classical part of $\taufed$
and let $q_0 \in Q$. If $q^1, \ldots, q^n$ are normal (geodesic)
coordinates around $q_0$ with respect to $\nabla_0$ then for any
$\alpha_{q_0} \in T^*_{q_0}Q$ the induced bundle coordinates
are normal Darboux coordinates
(see e.~g.~\cite[Sec. 2.5, p. 67] {Fed96} for definition)
around $\alpha_{q_0}$ with respect to $\nabla^0$ and
\BEQ {ClassFT}
    \left. \taufedcl (f) \right|_{\alpha_{q_0}}
    =
    \sum_{r=0}^\infty \frac{1}{r!}
    \left.\frac{\partial^r f} {\partial x^{i_1} \cdots \partial x^{i_r}}
    \right|_{\alpha_{q_0}}
    dx^{i_1} \vee \cdots \vee dx^{i_r}
\EEQ
for any $f \in C^\infty (T^*Q)$ where
$(x^1, \ldots, x^{2n}) = (q^1, \ldots, q^n, p_1, \ldots, p_n)$.
\end {THEOREM}
\begin {PROOF}
The fact that the bundle coordinates are normal Darboux
coordinates is proved in lemma \ref {NormalDarbouxLem}.
First we notice that it is sufficient to
prove (\ref {ClassFT}) for functions polynomial in the momenta
only. For those functions we shall prove the theorem by induction
on the order $k$ in the momenta. For $k=0$ the statement is true due
to corollary \ref {TauWeylPiCor} and theorem \ref {ClassFedTaylorTheo}
since the local expression of $\pi^* \tau_0 (\cdot)$ in normal (Darboux)
coordinates are clearly given by (\ref {ClassFT}) due to
lemma \ref {TaylorLem}.
Hence let $f \in C^\infty_{pp, k} (T^*Q)$. We prove the theorem
by a local argument: For the
coordinate function $\pi^*q^i$ we have
$\taufedcl (\pi^* q^i)|_{\alpha_{q_0}} =
(q^i + dq^i)|_{\alpha_{q_0}} = dq^i$ which
implies by (\ref {DweylClassic}) and the fact that $q^i$ is a
Hamiltonian for the Hamiltonian vector field $-\partial_{p_i}$ that
$\taufedcl (\frac{\partial f}{\partial p_i})|_{\alpha_{q_0}} =
i_s (\partial_{p_i}) \taufedcl (f)|_{\alpha_{q_0}}$. On the other
hand we compute $\Lie_{\partial_{p_i}} \taufedcl (f)$ at $\alpha_{q_0}$
using (\ref {iaXDweyl}) and (\ref {DweylClassic}) which implies
$(\Lie_{\partial_{p_i}} \taufedcl (f))|_{\alpha_{q_0}}
= i_s (\partial_{p_i}) \taufedcl (f)|_{\alpha_{q_0}}$ since
$\nabla^0_{X} \partial_{p_i} |_{\alpha_{q_0}} = 0$ for all
$X \in \Gamma (T(T^*Q))$. The homogeneity of $\taufed$ implies
$\Lie_\xi \taufedcl = \taufedcl \Lie_\xi$ and hence
$\taufedcl (f)$ is a polynomial in $p_i$ and $dp_i$
of maximal degree $k$. Then the last consideration implies
that at $\alpha_{q_0}$ it only depends on the combination $p_i + dp_i$.
Now using (\ref {ClassFT}) for
$\frac{\partial f}{\partial p_i} \in C^\infty_{pp, k-1} (T^*Q)$
due to lemma \ref {WichtigLem} the induction is easily finished.
\end {PROOF}

\section {Homogeneous Fedosov star products of standard ordered
          type}
\label {StandardSec}

In this section we shall construct a Fedosov star product of
standard ordered type. From now on we shall only use the
connection $\nabla^0$ on $T^*Q$ corresponding to a connection
$\nabla_0$ on $Q$ as given in (\ref {LiftedConnDef}). We define
the fibrewise standard ordered product for $a, b \in \WL$ by
\BEQ {FibStdDef}
    a \std b := \mu \circ e^{\frac{\lambda}{\im} i_s (\partial_{p_k})
                \otimes i_s (\partial^h_{q^k})} a \otimes b
\EEQ
where $\partial^h_{q^k} = \partial_{q^k}
+ (\pi^*\Gamma^l_{kr}) p_l \partial_{p_r}$ is the horizontal lift
of $\partial_{q^k}$ viewed as vector field on $Q$ with respect to
$\nabla_0$ to a vector field on $T^*Q$
(see (\ref {LiftVectDef})).
Clearly $\std$ is globally defined and an
associative deformation of $\mu$ which is of the form
(\ref{FibProd}). Moreover we define the maps
\BEQ {DeltaSDef}
    \Deltafib := i_s (\partial_{p_k}) i_s (\partial^h_{q^k})
    \quad
    \mbox { and }
    \quad
    \mathcal S := e^{\frac{\lambda}{2\im} \Deltafib}
\EEQ
which are again globally defined. Then the well-known
equivalence (see e.~g.~\cite {AW70}) of the standard ordered product
and the Weyl product in $\mathbb R^{2n}$ can easily be transfered
to a fibrewise equivalence: For $a, b \in \WL$ we have
\BEQ {FibWeylEquiStd}
    \mathcal S (a \fed b) = (\mathcal Sa) \std (\mathcal Sb)
\EEQ
and moreover, $\mathcal S$ commutes with $\mathcal H$ and $\Deg$
and thus $\mathcal H$ is a derivation of $\std$ and
$\std$ is homogeneous, too.
For the supercommutators using $\std$ we have
$\ads (a) = \mathcal S \circ \adfed (\mathcal S^{-1}a)
\circ \mathcal S^{-1}$
and $\ads (a) = 0$ iff $\adfed (a) = 0$ iff $\degs a = 0$.
Moreover we have
$[\delta, \Deltafib] = [\delta, \mathcal S] = [\delta, \mathcal S^{-1}] = 0$.
By a direct calculation of the commutator
$\mathcal B:= \frac{\im\lambda}{2} [\nabla^0, \Deltafib]$ we
obtain the following local expression for $\mathcal B$
\BEQ {BDefNablaDeltaCom}
    \mathcal B = (1\otimes dq^i) \frac{\im\lambda}{3} p_l
    \left(\pi^* R^l_{jik}\right) 
    i_s (\partial_{p_j}) i_s (\partial_{p_k})
\EEQ
where $R^l_{jik}$ are the components of the curvature tensor
of $\nabla_0$.
Now we conjugate the $\fed$-super\-de\-ri\-va\-tion $\nabla^0$
with $\mathcal S$ to obtain a $\std$-superderivation:
\BEQ {SnablaSinv}
    \mathcal S \nabla^0 \mathcal S^{-1} = \nabla^0 + \mathcal B
\EEQ
Note that no higher terms occur since $\mathcal B$ already commutes
with $\Deltafib$. Now $\nabla^0 + \mathcal B$ is a superderivation of $\std$
of antisymmetric degree $1$ and total degree $0$ which commutes
again with $\mathcal H$.

Now we shall follow two ways to obtain a Fedosov derivation for
$\std$: Firstly we just conjugate $\Dfed$ by $\mathcal S$ which will lead
indeed to a Fedosov derivation for $\std$ and secondly we start the
recursion new with $\nabla^0 + \mathcal B$ as term of total degree $0$ as
in theorem \ref {GenFDerivTheo}. The following proposition is
straight forward:
\begin {PROPOSITION} \label {DummyDProp}
The map $\mathcal D': \WL \to \WL$ defined by
\BEQ {DummyDDef}
    \mathcal D' := \mathcal S \Dfed \mathcal S^{-1}
\EEQ
is a superderivation of antisymmetric degree $1$ of $\std$ which
has square zero: ${\mathcal D'}^2 = 0$. It commutes with
$\mathcal H$ and we have
\BEQ {DummyD}
    \mathcal D' = - \delta + \nabla^0 + \mathcal B +
                  \frac{\lambda}{\im} \ads (\mathcal S\rfed)
\EEQ
where $\mathcal S\rfed = \rfed + \frac{\lambda}{2\im} \Deltafib \rfed$.
Hence $\mathcal D'$ satisfies the conditions of theorem
\ref {GenFTaylorTheo} and the corresponding Fedosov-Taylor series
$\tau'$ commutes with $\mathcal H$, too. Moreover for
$f, g \in C^\infty (T^*M)[[\lambda]]$
\BEQ {DummyStarDef}
    f *' g := \sigma \left(\tau' (f) \std \tau' (g)\right)
\EEQ
defines a homogeneous star product for $T^*Q$.
\end {PROPOSITION}
Moreover $\tau'$ satisfies a recursion formula analogously to the 
recursion formula (\ref {FedTaylorRecurs}) for $\taufed$ with the 
modification that $\nabla^0$ is replaced by $\nabla^0 + \mathcal B$, 
$\adfed$ by $\ads$ and $\rfed$ by $\mathcal S\rfed$.

Now the star product $*'$ is no longer of Weyl type but we have
the following analogue to the usual standard ordered product in
$\mathbb R^{2n}$:
\begin {DEFINITION}
A star product $*$ for a cotangent bundle $\pi: T^*Q \to Q$ is called
of standard ordered type iff for all
$\chi \in C^\infty (Q)[[\lambda]]$ and all
$f \in C^\infty (T^*Q)[[\lambda]]$
\BEQ {StandardOrderedTypeDef}
    (\pi^*\chi) * f = (\pi^* \chi) f .
\EEQ
\end {DEFINITION}
Then the standard ordered star products can be characterized by
the following easy proposition:
\begin {PROPOSITION}
Let $*$ be a star product for a cotangent bundle
$\pi: T^*Q \to Q$ written as
$f * g = \sum_{r=0}^\infty \lambda^r M_r (f, g)$
with bidifferential operators $M_r$ then the following statements
are equivalent:
\begin {enumerate}
\item $*$ is of standard ordered type.
\item For all $f, g \in C^\infty (T^*Q)[[\lambda]]$ and
      $\chi \in C^\infty (Q)[[\lambda]]$ we have
      $((\pi^* \chi) f) * g = (\pi^* \chi) (f*g)$.
\item In any local bundle chart the bidifferential operators
      $M_r$, $r \ge 1$ are of the form
      \BEQ {StandardMr}
          M_r (f, g) =
          \sum_{l, s, t}
          M^{j_1 \ldots j_l}_{i_1 \ldots i_s k_1, \ldots k_t}
          \frac{\partial^s f}
          {\partial p_{i_1} \cdots \partial p_{i_s}}
          \frac{\partial^{l+t} g}
          {\partial q^{j_1} \cdots \partial q^{j_l}
          \partial p_{k_1} \cdots \partial p_{k_t}} .
      \EEQ
\end {enumerate}
\end {PROPOSITION}
The following corollary is immediately checked using the obvious
fact that $\tau' \pi^* = \pi^* i^* \tau' \pi^*$ and the particular
form of $\std$ and lemma \ref {WichtigLem}:
\begin {COROLLARY} \label {DummyStarStandardCor}
On every cotangent bundle there exists a homogeneous
star product of standard ordered type namely the homogeneous
Fedosov star product $*'$.
\end {COROLLARY}

Now we shall discuss the other important alternative: starting
the recursion new with $\nabla^0 + \mathcal B$ as term of total degree $0$.
First we have
to show that $\nabla^0 + \mathcal B$ satisfies indeed the conditions of
theorem \ref {GenFDerivTheo}. Clearly $[\delta, \nabla^0 + \mathcal B] = 0$
and moreover
\[
    \frac{1}{2} \left[ \nabla^0 + \mathcal B, \nabla^0 + \mathcal B \right]
    = \mathcal S \nabla^0 \mathcal S^{-1} \mathcal S \nabla^0 \mathcal S^{-1}
    = \mathcal S \frac{\im}{\lambda} \adfed (R) \mathcal S^{-1}
    = \frac{\im}{\lambda} \ads (\mathcal SR)
    = \frac{\im}{\lambda} \ads (R)
\]
since $\mathcal SR = R + \frac{\lambda}{2\im} \Deltafib R$ 
and $\degs \Deltafib R = 0$ which implies $\ads (\Deltafib R) = 0$. 
Hence $(\nabla^0 + \mathcal B)^2$ is
an inner superderivation with the element $R$. It remains to show
that $(\nabla^0 + \mathcal B)R = 0$ but this is clear since $\Lie_\xi R = R$
and thus $\mathcal BR = 0$ due to lemma \ref {WichtigLem} and
(\ref {BDefNablaDeltaCom}) and $\nabla^0 R = 0$ anyway. Thus we can apply
indeed theorem \ref {GenFDerivTheo} and obtain the following
proposition completely analogously to proposition
\ref {DummyDProp} and corollary \ref {DummyStarStandardCor}:
\begin {PROPOSITION} \label {DstdProp}
There exists a unique element $\rstd \in \WL^1$ such that
$\delta \rstd = R + (\nabla^0 + \mathcal B) \rstd + \frac{\im}{\lambda}
\rstd \std \rstd$ and $\delta^{-1} \rstd = 0$.
Moreover $\mathcal H \rstd = \rstd$. Then the corresponding
Fedosov derivation
\BEQ {DstdDef}
    \Dstd = -\delta + \nabla^0 + \mathcal B + \frac{\im}{\lambda} \ads
    (\rstd)
\EEQ
has square zero and $\Dstd$ as well as the corresponding Fedosov-Taylor 
series $\taustd$ commute with $\mathcal H$. Then
\BEQ {StdStarDef}
    f \stars g := \sigma \left( \taustd (f) \std \taustd (g)
    \right)
\EEQ
where $f, g \in C^\infty (T^*Q)[[\lambda]]$ defines
a homogeneous star product of standard ordered type.
\end {PROPOSITION}
Again we have a recursion formula for $\taustd (f)$ analogously 
to (\ref {FedTaylorRecurs}) resp. (\ref {GenTaylorRecurs}) with the 
obvious modifications.
We shall refer to $\stars$ as the homogeneous Fedosov star product
of standard ordered type. At last we consider the element $\rstd$
more closely and notice
first that $\mathcal B\rstd = 0$ which follows immediately from the local
expression (\ref {BDefNablaDeltaCom}) for $\mathcal B$ and
$\mathcal H \rstd = \rstd$ and lemma \ref {WichtigLem}.
This implies that the recursion formula for $\rstd$ as proposed by
theorem \ref {GenFDerivTheo} can be simplified to
\[
    \rstd^{(k+3)} = \delta^{-1} \left( \nabla^0 \rstd^{(k+2)}
    - \frac{1}{2} \sum_{l=1}^{k-1} \left\{
    \rstd^{(l+2)}, \rstd^{(k-l+2)} \right\}_{\rm fib} \right)
\]
since $\mathcal H \rstd = \rstd$ and using the explicit expression for
$\std$ and lemma \ref {WichtigLem}. Thus $\rstd$ satisfies the
same recursion formula as $\rfed$ with the same first term
$\rstd^{(3)} = \delta^{-1} R = \rfed^{(3)}$ which implies that
they coincide:
\begin {LEMMA}
Let $\rstd$ be given as in proposition \ref {DstdProp} then
\BEQ {rstdrweyl}
    \rstd = \rfed .
\EEQ
\end {LEMMA}

Finally we compute the Fedosov-Taylor series $\tau'$ and $\taustd$
of the pull-back of a function on $Q$:
\begin {PROPOSITION} \label {DummyDpiETCProp}
Let $\chi \in C^\infty (Q)[[\lambda]]$ and $\mathcal D'$, $\tau'$ and
$\Dstd$, $\taustd$ be given as above. Then
\BEQ {DstdPi}
    \mathcal D' \pi^* = \Dstd \pi^* = \pi^* \mathcal D_0
\EEQ
\BEQ {taustdpi}
    \tau' (\pi^* \chi) = \taustd (\pi^* \chi) = \pi^* \tau_0 (\chi)
\EEQ
where $\mathcal D_0$ as in (\ref {ClassFDerivDef}) and $\tau_0$
is the formal Taylor series with respect to $\nabla_0$ as in
(\ref {ClassTau}).
\end {PROPOSITION}
\begin {PROOF}
Clearly $\mathcal B\pi^* = 0$ and
$\ads (\mathcal S\rfed)$ as well as $\ads (\rfed)$ applied to pull-backs
by $\pi^*$ reduce to fibrewise Poisson brackets due to
lemma \ref {WichtigLem} since $\Lie_\xi \rfed = \rfed$.
This implies $\mathcal D' \pi^* = \Dstd \pi^* = \Dfed \pi^*$ and
by (\ref {DweylPiDNull}) the first equation is proved which
implies the second equation.
\end {PROOF}

\section {Equivalence of $\starf$, $*'$, and $\stars$}
\label {EquiSec}

In this section we shall construct equivalence transformations of
the star products $\starf$, $*'$ and $\stars$. Though DeWilde and
Lecomte had shown in \cite [Proposition 4.2.] {DL83a} that any
homogeneous star products of {\em Weyl type} on a cotangent bundle
are equivalent this is in so far a non-trivial problem since
the star products $*'$ and $\stars$ are evidently not of Weyl type.
First we define for
$f \in C^\infty (T^*Q)[[\lambda]]$
\BEQ {TDef}
    Tf := \sigma \left( \mathcal S^{-1} \tau' (f) \right)
\EEQ
then the fibrewise equivalence of $\fed$ and $\std$ induces via
$T$ an equivalence of $\starf$ and $*'$:
\begin {PROPOSITION} \label {TEquivProp}
The map $T$ is homogeneous $[\mathcal H, T] = 0$ and an
equivalence transformation between $*'$ and $\starf$ 
\BEQ {StarWEquiv}
    T (f *' g) = Tf \starf Tg
\EEQ
and we have $\taufed (Tf) = \mathcal S^{-1} \tau' (f)$ and
$T^{-1} f = \sigma \left(\mathcal S\taufed (f)\right)$
where $f, g \in C^\infty (T^*Q)[[\lambda]]$. Moreover $T$ can be
written as formal series
$T = \id + \sum_{r=1}^\infty \lambda^r T_r$ where $T_r$ is a
differential operator of order $2r$ given by
\BEQ {Tr}
    T_r f = \sum_{s=0}^r \frac{1}{s!} \left(\frac{\im}{2}\right)^s
            \sigma \left(\Deltafib^s \tau' (f)^{(2r)}_{2s} \right)
\EEQ
where $\tau' (f)^{(k)}_l$ denotes the term of total degree $k$ and
symmetric degree $l$.
\end {PROPOSITION}
\begin {PROOF}
The first part follows from proposition \ref {EquivalenceTrafoProp}
and the order of differentiation follows form
the fact that $\tau'(f)^{(2r)}$ is a differential operator of
order $2r$ which follows directly from the recursion formula 
(\ref {GenTaylorRecurs}) for $\tau'$.
\end {PROOF}

The construction of an equivalence transformation between $*'$ and
$\stars$ is more involved. First we prove that the Fedosov
derivations $\mathcal D'$ and $\Dstd$ are conjugated to each other
by a fibrewise automorphism of $\std$:
\begin {THEOREM} \label {hAutoTheo}
Let $\Dstd$ and $\mathcal D'$ be the Fedosov derivations constructed
as in section \ref {StandardSec} then there exists an element
$h \in \W(Q)_{\rm cl}$ such that
\BEQ {DstdConjugate}
    \Dstd = e^{\ads (\pi^* h)} \mathcal D' e^{-\ads (\pi^* h)}
          = e^{\ads (\pi^* h)} \mathcal S \Dfed \mathcal S^{-1}
            e^{-\ads (\pi^*h)}.
\EEQ
\end {THEOREM}
\begin {PROOF}
Let $h \in \W(Q)_{\rm cl}$ be an arbitrary element then by direct
calculation using the properties of $\std$ and lemma \ref {WichtigLem}
we obtain
\[
    e^{\ads (\pi^* h)} (-\delta) e^{-\ads (\pi^* h)}
    =
    -\delta + \ads \left(\pi^* (\delta_0 h\right))
\]
\[
    e^{\ads (\pi^* h)} (\nabla^0 + \mathcal B) e^{-\ads (\pi^*  h)}
    =
    \nabla^0 + \mathcal B - \ads \left(\pi^* \nabla_0 h\right)
\]
\[
    e^{\ads (\pi^* h)} \frac{\im}{\lambda} \ads (\mathcal S\rfed)
    e^{-\ads (\pi^* h)}
    =
    \frac{\im}{\lambda} \ads\left(
    \rfed - \frac{\lambda}{\im} \pi^*\left(i_s (\varrho_0) h\right)
    + \frac{\lambda}{2\im} \pi^* \left (\tr \varrho_0\right)
    \right)
\]
where $\tr := i_s (\partial_{q^k}) i(dq^i)$ and
$\varrho_0$ is given as in theorem \ref {ClassFedDerivTheo}.
Collecting these results we obtain due to $\rstd = \rfed$
\[
    e^{\ads (\pi^* h)} \mathcal D' e^{-\ads (\pi^* h)}
    =
    \Dstd - \ads \left( \pi^*\left( \mathcal D_0 h - \frac{1}{2}
    \tr \varrho_0\right) \right)
\]
where $\mathcal D_0$ as in theorem \ref {ClassFedDerivTheo} applied
for $\nabla_0$ on $Q$. Hence we have to find an element $h$ such
that the second term vanishes which is the case iff we find a
one-form $\alpha \in \Gamma (\bigwedge^1 T^*M)$ such that
\BEQ {DNullh}
    \mathcal D_0 h = \frac{1}{2}
    \left(\tr\varrho_0 + 1 \otimes \alpha \right)
\EEQ
since in this case $\mathcal D_0 h - \frac{1}{2} \tr \varrho_0$ is
a central element. Note that we only consider
elements $h \in \W_{\rm cl}$.
A necessary condition for (\ref {DNullh}) to
have a solution is that the right hand side is
$\mathcal D_0$-closed since ${\mathcal D_0}^2 = 0$.
We have $\mathcal D_0 \varrho_0 = R_0$ due to theorem
\ref {ClassFedDerivTheo}
and apply $\tr$ on both sides leading to
$\tr \mathcal D_0 \varrho_0 = \tr R_0$
where $R_0$ is the curvature tensor of $\nabla_0$. Now a straight forward
computation leads to
$\tr \mathcal D_0 \varrho_0 = \mathcal D_0 \tr \varrho_0$.
On the other hand if $\alpha \in \Gamma (\bigwedge^1 T^*Q)$ then
$\mathcal D_0 (1 \otimes \alpha) = 1 \otimes d \alpha$ since
$\nabla_0$ is torsion-free. Since the trace of the curvature tensor
$R_0$ is an exact two-from (see lemma \ref {TrRExactLem}) we find
always a one-form $\alpha$ such that $\tr R_0 = - d \alpha$ and hence
$\mathcal D_0 (\tr \varrho_0 + 1 \otimes \alpha) = 0$ iff
$\alpha$ satisfies $\tr R_0 = - d\alpha$ and thus the necessary
condition is fulfilled. But this is also sufficient
since the $\mathcal D_0$-cohomology is trivial on forms:
Indeed we define $\WL^+ := \{ a \in \WL | \sigma (a) = 0\}$ then we have
the following lemma (see \cite [Theorem 5.2.5] {Fed96}):
\begin {LEMMA} \label {DNullCohomLem}
Let $\mathcal D_0^{-1} : \WL(Q) \to \WL(Q)$
be defined by
\BEQ {DNullInvDef}
    \mathcal D_0^{-1} := - \delta_0^{-1} \frac{1}{1 -
    [\delta_0^{-1}, \nabla_0 + i_s (\varrho_0)]}
\EEQ
then for any $a \in \WL^+(Q)$ the following
`deformed Hodge decomposition' holds
\BEQ {DefHodgeDecomp}
    \mathcal D_0 \mathcal D_0^{-1} a
    + \mathcal D_0^{-1} \mathcal D_0 a
    = a
\EEQ
and $[\delta_0^{-1}, \nabla_0 + i_s (\varrho_0)]$ commutes with
$\delta_0^{-1}$ and $\mathcal D_0$ and
$\cc {\mathcal D_0^{-1} a} = \mathcal D_0^{-1} \cc a$.
\end {LEMMA}
Note that $\mathcal D_0^{-1}$ is a well-defined formal series
in the symmetric degree. Now choose a one-form $\alpha$ with
$\tr R_0 = - d\alpha$. Then we define
\BEQ {hDef}
    h := \frac{1}{2} \mathcal D_0^{-1}
    \left(\tr \varrho_0 + 1 \otimes \alpha \right)
\EEQ
and notice that $\sigma (h) = 0$ since $\mathcal D_0^{-1}$ raises
the symmetric degree. Moreover 
$\sigma\left(\tr \varrho_0 + 1\otimes\alpha \right)=0$ 
hence we can apply the lemma and find
$\mathcal D_0 h = \frac{1}{2} \left(\tr \varrho_0 + 1 \otimes \alpha\right)$
and thus the theorem is proved.
\end {PROOF}
\begin {COROLLARY}
If $h$ is a solution of (\ref {DNullh}) for a fixed one-form
$\alpha$ then every other solution $h'$ is obtained by
$h' = h + \tau_0 (\varphi)$ with $\varphi \in C^\infty (Q)$. For a
fixed one-form $\alpha$ satisfying $\tr R_0 = - d\alpha$
there exists a unique solution of (\ref {DNullh}) with
$\sigma (h) = \varphi$ for any $\varphi \in C^\infty (Q)$ namely
$h = \frac{1}{2} \mathcal D_0^{-1}
(\tr \varrho_0 + 1 \otimes \alpha) + \tau_0 (\varphi)$.
The one-from $\alpha$ is determined up to a closed one-form and can be
chosen to be real $\cc \alpha = \alpha$ which leads to a real
$h = \cc h$ iff $\varphi$ is real.
\end {COROLLARY}
\begin {COROLLARY}
If the connection $\nabla_0$ is unimodular then there exists a
canonical solution $h$ of (\ref {DstdConjugate}) uniquely
determined by $\mathcal D_0 h = \frac{1}{2}\tr \varrho_0$ and
$\sigma (h) = 0$ namely
$h = \frac{1}{2} \mathcal D_0^{-1} \tr \varrho_0$.
In this case $h = \cc h$ is real.
\end {COROLLARY}
\begin {COROLLARY}
Let $h$ be an arbitrary solution of (\ref {DNullh}) then
$e^{\ads(\pi^* h)}$ is a fibrewise automorphism of $\std$
satisfying (\ref {DstdConjugate}) and commuting with $\mathcal H$.
\end {COROLLARY}
Now we can use the fibrewise automorphism $e^{\ads (\pi^* h)}$
to construct an equivalence transformation between $\stars$ and $*'$ 
analogously to the construction of $T$ in (\ref {TDef}).
We define for $f \in C^\infty (T^*Q)[[\lambda]]$
\BEQ {VDef}
    Vf := \sigma \left( e^{-\ads (\pi^* h)} \taustd (f) \right)
\EEQ
and get the following proposition completely analogously to
proposition \ref {TEquivProp}:
\begin {PROPOSITION}
For any solution $h$ of (\ref {DstdConjugate}) the map $V$ is
homogeneous $[\mathcal H, V] = 0$ and an
equi\-valence transformation between $\stars$ and $*'$:
\BEQ {StarSEquiv}
    V (f \stars g) = Vf *' Vg
\EEQ
and we have $\tau' (Vf) =  e^{-\ads (\pi^* h)} \taustd (f)$ and
$V^{-1} f = \sigma \left(e^{\ad(\pi^* h)} \tau' (f)\right)$
where $f, g \in C^\infty (T^*Q)[[\lambda]]$.
\end {PROPOSITION}
To compute the orders of differentiation in $V$ we first need the
following lemma obtaining by the way that $\stars$ is a Vey product:
\begin {LEMMA}
Let $\taustd$ be the Fedosov-Taylor series as in proposition
\ref {DstdProp} and let $\taustd (\cdot)^{(k)}_l$ be the term of
total degree $k$ and symmetric degree $l$. Then
$\taustd (\cdot)^{(2r)}_{2s}$ resp. $\taustd (\cdot)^{(2r+1)}_{2s+1}$
is a differential operator of order $r+s$ resp. $r+s+1$ for all
$r, s \in \mathbb N$. Moreover this implies that the homogeneous
Fedosov star product of standard ordered type $\stars$ is a Vey product.
\end {LEMMA}
\begin {PROOF}
This lemma is proved by a straight forward induction on the total
degree using the recursion formula (\ref {GenTaylorRecurs}) for
$\taustd$. Then the Vey type property of $\stars$ follows
directly.
\end {PROOF}
\begin {COROLLARY}
For any solution $h$ of (\ref {DNullh}) the equivalence
transformation $V$ can be written as formal series
$V = \id + \sum_{r=1}^\infty \lambda^r V_r$ where $V_r$ is a
differential operator of order $r$.
\end {COROLLARY}

\section {The standard representation}
\label {RepSec}

Now we shall construct a canonical representation of the
fibrewise algebra $(\W, \std)$ and of the star product algebra
$(C^\infty (T^*Q)[[\lambda]], \stars)$ reproducing the well-known
standard order quantization rule for cotangent bundles.
First we define the representation space
\BEQ {FibHDef}
    \mathfrak H := \W(Q)
\EEQ
and define the {\em fibrewise standard ordered represenation}
for $a \in \W$ and $\Psi \in \mathfrak H$ by
\BEQ {FibStdRepDef}
    \fibsrep (a) \Psi := i^* \left( a \std \pi^* \Psi \right)
\EEQ
and notice that $\fibsrep: \W \to \End (\mathfrak H)$ is indeed
a representation of $\W$ with repect to $\std$:
\begin {LEMMA}
Let $\std$ be the fibrewise standard ordered product and
$\fibsrep$ be defined as in (\ref {FibStdRepDef}) then $\fibsrep$
is a $\std$-representation of $\W$ on $\mathfrak H$, i.~e.~for
$a, b \in \W$
\BEQ {FibStdRep}
    \fibsrep (a \std b) = \fibsrep (a) \fibsrep (b) .
\EEQ
Furthermore
\BEQ {FibWeylRepDef}
    \fibfrep (a) := \fibsrep (\mathcal Sa)
\EEQ
defines a representation with respect to the fibrewise Weyl
product $\fed$ of $\W$ on $\mathfrak H$ given by
\BEQ {FibWeylRep}
    \fibfrep (a)\Psi = i^* \mathcal S\left( a \fed \pi^* \Psi \right) .
\EEQ
\end {LEMMA}
\begin {PROOF}
The $\mathbb C[[\lambda]]$-linearity of $\fibsrep$ is obvious and the
representation property is proved by straight forward computation.
Then the fibrewise equivalence (\ref {FibWeylEquiStd}) of $\std$
and $\fed$ implies that $\fibfrep$ is a representation
with respect to $\fed$.
\end {PROOF}

Now we shall construct a representation $\srep$ for the star product
$\stars$ induced by $\fibsrep$. First we notice that the restriction
of $\fibsrep$ to $\ker\Dstd \cap \W$ is still a
representation of $\ker\Dstd \cap \W$ and hence it induces via
$\taustd$ a representation of $C^\infty (T^*Q)[[\lambda]]$ on
$\mathfrak H$. But we have in mind to construct a representation of
$\ker\Dstd \cap \W$ on the smaller representation space
$i^*(\ker\Dstd \cap \pi^*\mathfrak H) \subset \mathfrak H$ which is via 
$\tau_0$ in bijection with $C^\infty (Q)[[\lambda]]$. In the case of the
standard ordered product $\stars$ this can be done directly:
\begin {THEOREM}
Let $\Dstd$ be the Fedosov derivation constructed as in
proposition \ref {DstdProp} and let $\taustd$ be the corresponding
Fedosov-Taylor series and $\stars$ the homogeneous Fedosov star
product of standard ordered type. Then
\BEQ {DstdpiiCom}
    \Dstd \pi^* i^* = \pi^* i^* \Dstd
\EEQ
and
\BEQ {StdRepDef}
    \srep (f) \psi := i^* (f \stars \pi^* \psi) =
    \sigma (\fibsrep (\taustd (f)) \tau_0 (\psi))
\EEQ
where $f \in C^\infty (T^*Q)[[\lambda]]$ and
$\psi \in C^\infty (Q)[[\lambda]]$ defines a representation of
$C^\infty (T^*Q)[[\lambda]]$ with respect to $\stars$
on $C^\infty (Q)[[\lambda]]$.
\end {THEOREM}
\begin {PROOF}
One immediately checks that $\delta$, $\nabla^0$ and $\mathcal B$
commute with $\pi^* i^*$. Moreover $\ads (\rstd)$ commutes
with $\pi^* i^*$ too due to the particular form of $\std$ and
$\Lie_\xi \rstd = \rstd$ and lemma \ref {WichtigLem}.
Hence (\ref {DstdpiiCom}) is shown.
This ensures that $\srep$ defines indeed a representation which
is now a straight forward computation using (\ref {FibStdRep})
and proposition \ref {DummyDpiETCProp}.
\end {PROOF}

Note that (\ref {DstdpiiCom}) is crucial for the representation
property of $\srep$ and that neither $\Dfed$ nor $\mathcal D'$
commute with $\pi^* i^*$. We shall refer to $\srep$ as
{\em standard representation}
with respect to $\stars$. A representation with respect to the
homogeneous Fedosov star product of Weyl type $\starf$ can be
constructed by use of the equivalence transformations $V$ and $T$
but we shall see in the next section that there is another
star product $\starw$ of Weyl type with a corresponding
representation.
Now we compute an explicit formula for the representation $\srep$
and rediscover the well-known standard order quantization rule for
cotangent bundles (see (\ref {StandOrdCov})):
\begin {THEOREM} \label {StandardRepTheo}
Let $\stars$ be the homogeneous Fedosov star product of standard
ordered type and let $\srep$ be the corresponding standard
representation. Then for $f \in C^\infty (T^*Q)[[\lambda]]$ and
$\psi \in C^\infty (Q)[[\lambda]]$ we have
\BEQ {StdRep}
    \srep (f) \psi =
    \sum_{r=0}^\infty \frac{1}{r!}
    \left(\frac{\lambda}{\im}\right)^r
    i^*\left( \frac{\partial^r f}{\partial p_{i_1} \cdots
    \partial p_{i_r}} \right)
    i_s (\partial_{q^{i_1}}) \cdots i_s (\partial_{q^{i_r}})
    D_0^{(r)} \psi
\EEQ
where $D_0^{(r)}\psi 
:= \frac{1}{r!} (dq^k \vee {\nabla_0}_{\partial q^k})^r \psi$ 
is the $r$th symmetrized covariant derivative of $\psi$ with 
respect to $\nabla_0$. In particular for a function
$f \in C^\infty_{pp, k} (T^*Q)$ polynomial in the momenta of order
$k$ we have
\BEQ {StdRepPoly}
    \srep (f) \psi = \frac{1}{k!} \left(\frac{\lambda}{\im}\right)^k
    \left< F, D_0^{(r)}\psi \right>
\EEQ
where $F \in \Gamma (\bigvee^k TQ)$ is the symmetric tensor field
such that $f = \widehat F$ and
$\left<\cdot, \cdot\right>$ denotes the natural pairing. If $\widehat X$ 
is linear in the momentum where $X \in \Gamma (TQ)$ then
\BEQ {StdRepVect}
    \srep (f) \psi = \frac{\lambda}{\im} \Lie_X \psi .
\EEQ
\end {THEOREM}
\begin {PROOF}
We compute $\srep (f) \psi$ using (\ref {StdRepDef}) and
proposition \ref{DummyDpiETCProp} and obtain
\[
    \srep (f) \psi = \sum_{r=0}^\infty \frac{1}{r!}
                     \left( \frac{\lambda}{\im} \right)^r
                     i^* \sigma \left(
                     i_s (\partial_{p_{i_1}}) \cdots
                     i_s (\partial_{p_{i_r}}) \taustd (f)
                     \right)
                     \sigma\left(
                     i_s (\partial_{q^{i_1}}) \cdots
                     i_S (\partial_{q^{i_r}}) \tau_0 (\psi)
                     \right) .
\]
But since $\tau_0$ is the Taylor series with respect to $\nabla_0$
we find that the last part
$\sigma( i_s (\partial_{q^{i_1}}) \cdots i_s
(\partial_{q^{i_r}}) \tau_0 (\psi))$
equals
$i_s (\partial_{q^{i_1}}) \cdots i_s (\partial_{q^{i_r}}) D^{(r)}_0 \psi$
and thus we only have to prove that
$i^* \sigma (i_s (\partial_{p_{i_1}}) \cdots
i_s (\partial_{p_{i_r}}) \taustd (f))$ coincides with
$i^* (\frac{\partial^r f}{\partial p_{i_1} \cdots
\partial_{p_{i_r}}})$. But this follows directly from the
following lemma:
\end {PROOF}
\begin {LEMMA}
Let $f \in C^\infty (T^*Q)$ then in any local bundle chart there
exist locally defined elements $\tilde \tau_i^r (f)$ such that
for any total degree $r \ge 0$
\BEQ {tildetau}
    \taustd (f)^{(r)} = D^{(r)} f + dq^i \vee \tilde\tau_i^r (f)
\EEQ
\end {LEMMA}
\begin {PROOF}
This lemma is proved by a straight forward induction on the total
degree using the recursion formulas (\ref {GenTaylorRecurs}) for $\taustd$.
\end {PROOF}
\begin {COROLLARY} \label {SRepInjCor}
The restriction of $\srep$ to the $\mathbb C[\lambda]$-submodule
$C^\infty_{pp} (T^*Q)[\lambda]$ as well as the
restriction to the $\mathbb C[[\lambda]]$-submodule
$C^\omega_{p} (T^*Q)[[\lambda]]$ is injective.
\end {COROLLARY}

\section {Two different homogeneous star products of Weyl type}
\label {WeylRepSec}

In this section we shall construct an analogue to the operator $N$
mentioned in the introduction. This operator allows us to define
a star product of Weyl type $\starw$ equivalent to $\stars$ which
corresponds to the Weyl ordering prescription of flat 
$\mathbb R^{2n}$ and turns out to be different form $\starf$ in general.

In this section we sometimes denote the complex conjugation by
$\mathcal C$, i. e. $\mathcal C a := \cc a$ for $a \in \WL$ and
define
\BEQ {CCsDef}
    \CCs := e^{\ads (\pi^* h)} \mathcal {SCS}^{-1} e^{-\ads (\pi^*h)}
\EEQ
where $h \in \W (Q)_{\rm cl}$ is constructed as in theorem
\ref {hAutoTheo} with $\sigma (h) = 0$ incorporating a particular but
fixed choice of a {\em real} one-form
$\alpha \in \Gamma (\bigwedge^1 T^*Q)$ such that
$-d\alpha = \tr R_0$. Then we have the following lemmata:
\begin {LEMMA}
Let $a, b \in \WL$ with $\dega a = ka$ and $\dega b = lb$ then
\BEQ {CAntiAuto}
    \mathcal C (a \fed b) = (-1)^{kl} (\mathcal C b) \fed (\mathcal C a)
\EEQ
and $[\Dfed, \mathcal C] = 0$ and
\BEQ {CCsAntiAuto}
    \CCs (a \std b) = (-1)^{kl} (\CCs b) \std (\CCs a)
\EEQ
\BEQ {CCsDsCCsHom}
    \left[ \Dstd, \CCs\right] = 0 = \left[\mathcal H, \CCs \right] .
\EEQ
\end {LEMMA}
\begin {PROOF}
Equation (\ref {CAntiAuto}) and the fact $[\Dfed, \mathcal C] = 0$ are
well-known (see e.~g.~\cite [Lemma 3.3] {BW96a}) and imply together with
(\ref {DstdConjugate}) the other statements of the lemma.
\end {PROOF}
\begin {LEMMA} \label{CSantiaut}
Let $f, g \in C^\infty (T^*Q)[[\lambda]]$ then the map
\BEQ {CsDef}
    \Cs f := \sigma \left( \CCs \taustd (f) \right)
\EEQ
defines an involutive anti-$\mathbb C[[\lambda]]$-linear 
anti-automorphism of the star product $\stars$, i.~e. $\Cs^2 = \id$ and
\BEQ {CsAntiAuto}
    \Cs (f \stars g) = (\Cs g) \stars (\Cs f)
\EEQ
and we have $\taustd (\Cs f) = \CCs \taustd (f)$ and
$[\Cs, \mathcal H] = 0$. Moreover $\Cs \mathcal C$ is a formal series of
differential operators
\BEQ {CsCSeries}
    \Cs \mathcal C = \id + \sum_{r=1}^\infty \lambda^r \Cs^{(r)}
\EEQ
where $\Cs^{(r)}$ is a differential operator of order $2r$ and
$\Cs^{(1)} = \frac{1}{\im} \Delta$ and
\BEQ {CovDeltaDef}
    \Delta := {\nabla^0}^2_{(\partial_{p_i}, \partial_{q^i}^h)} +
              \nabla^0_{(\alpha^v)}
\EEQ
which is clearly globally defined and takes the following form in a
bundle chart
\BEQ {CovDeltaLocal}
    \Delta = \partial_{q^i} \partial_{p_i}
             + p_r \pi^*(\Gamma^r_{ij}) \partial_{p_i} \partial_{p_j}
             + \pi^*(\Gamma^i_{ij}) \partial_{p_j}
             + \pi^*(\alpha_j) \partial_{p_j} .
\EEQ
\end {LEMMA}
\begin {PROOF}
Equation (\ref {CsAntiAuto}) is proved in a completely analogous
manner as proposition \ref {EquivalenceTrafoProp}. Since the term
$\taustd (\cdot)^{(k)}$ of total degree $k$ is easily seen to be of
order $k$ as a differential operator and since $\CCs$ does obviously
not decrease the total degree ($\mathcal S$ being of total
degree $0$ and $\pi^* h$ contains only terms of positive total degree)
equation (\ref {CsCSeries}) follows. Moreover $\Cs^{(1)} = -\im \Delta$
follows by straight forward computation using
$\taustd (f)^{(2)} 
= \frac{1}{2} (dx^k \vee {\nabla^0}_{\partial_{x^k}})^2 f$ 
due to (\ref {GenTaylorRecurs}) and $h_1 = - \frac{1}{2} \alpha \otimes 1$
where $h_1$ is the term of symmetric degree $1$ in $h$ due 
to (\ref {hDef}).
\end {PROOF}

Now we come to a simple analogue of the operator $N$ mentioned 
in (\ref {SFlat}):
\BEQ {NeumaierDef}
    N := \exp \left( \frac{\lambda}{2\im} \Delta\right)
\EEQ
which clearly is a formal power series in $\lambda$ of differential
operators and $[N, \mathcal H] = 0$. Motivated by equation
(\ref {FlatWeylStdEquiv}) we define an equivalent star product to
$\stars$ using $N$ as an equivalence transformation by
\BEQ {TrueWeylProdDef}
    f \starw g := N^{-1} \left( (Nf) \stars (Nf) \right)
\EEQ
for $f, g \in C^\infty (T^*Q)[[\lambda]]$ which is clearly
bidifferential and homogeneous. Furthermore we have the following
theorem:
\begin {THEOREM}
\begin {enumerate}
\item The operator $N^2 \mathcal C$ coincides with $\Cs$.
\item The star product $\starw$ is a star product of Weyl type
      where in particular
      \BEQ {CCTrueWeylProd}
          \cc{ f \starw g} = \cc g \starw \cc f
      \EEQ
      for $f, g \in C^\infty (T^*Q)[[\lambda]]$.
\item A representation $\wrep$ of $(C^\infty (T^*Q)[[\lambda]], \starw)$
      on $C^\infty (Q)[[\lambda]]$ is given by the following
      analogue of equation (\ref {WeylundStandard})
      \BEQ {WrepDef}
          \wrep (f) :=  \srep (Nf)  .
      \EEQ
\end {enumerate}
\end {THEOREM}
\begin {PROOF}
\begin {enumerate}
\item We proceed in two steps: Let us first prove that
      $N^2 \mathcal C$ is an involutive anti-linear anti-automorphism
      of $(C^\infty(T^*Q)[[\lambda]],\stars)$. This is equivalent to
      the identity
      \[
          N^2 \mathcal C(f\stars g) -
          (N^2 \mathcal C g)\stars (N^2 \mathcal C f) = 0
      \]
      for all $f, g \in C^\infty (T^*Q)[[\lambda]]$. Since $\stars$
      consists in bidifferential operators this identity is
      an identity of bidifferential operators for each power of $\lambda$.
      Hence it suffices to check it on an arbitrarily small
      neighbourhood of each point $m \in T^*Q$ on
      $f, g \in C_{pp}^\infty(T^*Q)$. Since the standard representation
      $\srep$ is injective on $C_{pp}^\infty(T^*Q)[\lambda]$
      it is sufficient to prove this identity after having
      applied $\srep$ to the left hand side. Now let $q \in Q$ be
      arbitrary and choose a contractible open neighbourhood $U$ around
      $q$ which lies in the domain of a chart. Suppose that the supports
      of $f, g$ lie both in $\pi^{-1}(U)$.
      Then there is a local volume form $\mu$ on $U$ such that
      ${\nabla_0}_X \mu = \alpha (X) \mu$ for each vector field $X$ on
      $U$: indeed, choose an arbitrary local volume form
      $\mu'$ on $U$. Then ${\nabla_0}_X \mu' = \alpha'(X) \mu'$ for a
      certain locally defined one-form $\alpha'$ on $U$. Since
      $d\alpha'(X, Y)\mu' =
      ([{\nabla_0}_X, {\nabla_0}_Y] - {\nabla_0}_{[X,Y]})\mu'
      = - (\tr R)(X, Y)\mu' = d\alpha(X, Y)\mu'$
      it follows from the Poincar\'e lemma that there is a local
      real-valued smooth function $\phi$ on $U$ such that
      $\alpha' = \alpha + d\phi$. Then $\mu := e^{-\phi} \mu'$ will
      clearly do the job. Consider next the space $\mathcal {D}(U)$ of
      all smooth complex-valued functions on $Q$ whose support lies
      in $U$ and is compact. This space is an inner
      product space with respect to the Lebesgue integral
      $\langle\phi, \psi\rangle_U := \int_U \cc\phi\psi \mu$. 
      We define the following covariant divergence operator
      $\Div_\alpha: \Gamma(\bigvee^k TQ) \to \Gamma(\bigvee^{k-1}TQ)$
      with $\Gamma(\bigvee^l TQ):=\{0\}$ for negative integers $l$:
      \BEQ {alphaDivDef}
          \Div_\alpha S := \Div S + i_s (\alpha) S 
          \qquad \mbox { where } \qquad
          \Div S := i_s(dq^i) \nabla_{\partial_{q^i}}S
      \EEQ
      which is clearly globally defined. Let $T$ be in
      $\Gamma(\bigvee^{k+1}TQ)$. Using the global coordinates
      $(q^1,\ldots,q^n)$ in $U$ and
      $D_0 := dq^i \vee {\nabla_0}_{\partial_i}$ we get for any
      $\psi\in\mathcal {D}(U)$:
      \begin {eqnarray*}
          \lefteqn{ \frac{1}{(k+1)!} T^{i_1\cdots i_{k+1}}
          i_s(\partial_{i_1}) \cdots i_s(\partial_{i_{k+1}})
          D_0^{k+1} \psi } \\
          & = & \Lie_{\partial_j} \left(
                \frac{1}{(k+1)!} T^{i_1\cdots i_{k+1}}
                i_s(\partial_{i_1})\cdots i_s(\partial_{i_{k+1}})
                \left(dq^j \vee D_0^k \psi \right) \right)  \\
          &   & - \frac{1}{(k+1)!}
                ({\nabla_0}_{\partial_j}T)^{i_1\cdots i_{k+1}}
                i_s(\partial_{i_1}) \cdots i_s(\partial_{i_{k+1}})
                \left(dq^j \vee D_0^k \psi \right)   \\
          &   & + \frac{1}{(k+1)!} \Gamma^r_{rj} T^{i_1\cdots i_{k+1}}
                i_s(\partial_{i_1}) \cdots i_s(\partial_{i_{k+1}})
                \left(dq^j \vee D_0^k \psi \right)  \\
          & = & \Lie_{\partial_j} \left(
                \frac{1}{k!} T^{ji_1\cdots i_{k}}
                i_s(\partial_{i_1})\cdots i_s(\partial_{i_{k}})
                D_0^k \psi \right)  \\
          &   & - \frac{1}{k!}
                ({\nabla_0}_{\partial_j}T)^{ji_1\cdots i_{k}}
                i_s(\partial_{i_1})\cdots i_s(\partial_{i_{k}})
                D_0^k \psi  \\
          &   & + \frac{1}{k!} \Gamma^r_{rj} T^{ji_1\cdots i_{k}}
                i_s(\partial_{i_1})\cdots i_s(\partial_{i_{k}})
                D_0^k \psi.
      \end {eqnarray*}
      Since for any vector field $X$ on $Q$ we obviously have
      $\Lie_X \mu = (\Div_\alpha X) \mu$ it can easily be seen by
      partial integration and induction that the following identity is
      true:
      \begin {eqnarray*}
          \lefteqn{
          \int_U \cc \phi \srep(\widehat{T})\psi \; \mu
          =
          \int_U \cc{\left(
          \frac{1}{k!} \left(\frac{\lambda}{\im}\right)^k
          (\Div_\alpha)^k (\phi \cc{T})\right)} 
          \; \psi \; \mu } \\
          & = & \int_U \cc{\left(
                \frac{1}{k!} \left(\frac{\lambda}{\im}\right)^k
                \sum_{s=0}^k {k \choose s}
                \left((\Div_\alpha)^s \cc{T} \right)^{i_1\cdots i_{k-s}}
                i_s(\partial_{i_1}) \cdots i_s(\partial_{i_{k-s}})
                {D_0}^{k-s} \phi\right)} \; \psi \; \mu \\
          & = & \int_U \cc{\left(
                \srep \left(\sum_{s=0}^k
                \frac{1}{s!} \left(\frac{\lambda}{\im}\right)
                \widehat{\left((\Div_\alpha)^s \cc{T}\right)} \right)\phi
                \right)} \; \psi \; \mu
                =
                \int_U \cc{\left(\srep\left(N^2\widehat{\cc T}\right)
                \phi \right)} \; \psi \; \mu
      \end{eqnarray*}
      where in the last equality we have used the identity
      $\widehat{(\Div_\alpha)^s T} = \Delta^s \widehat T$ for
      $s \in \mathbb N$ which follows from
      \BEQ {alpahDivIdentity}
          \widehat{\Div_\alpha T} = \Delta \widehat{T}
      \EEQ
      which can easily be proved by calculation. It follows that on
      $U$ the differential operator $\srep(\widehat{T})$
      has a formal adjoint $\srep(\widehat{T})^\dagger$ given by
      $\srep(N^2\widehat{\cc T})$. By the well-known rules of multiplying
      adjoints, viz.
      \[
         \left(\srep(\widehat{T})\srep(\widehat{S})\right)^\dagger
         = \srep(\widehat{S})^\dagger \srep(\widehat{T})^\dagger
      \]
      we finally get that $N^2 \mathcal C$ is an involutive anti-linear
      anti-automorphism since clearly
      $N^2 \mathcal C N^2 \mathcal C = \id$.

      Secondly, since now both $\Cs$ and $N^2 \mathcal C$ are involutive
      anti-linear anti-automorphisms the composition
      $\Cs^{-1} N^2 \mathcal C$ will be a linear automorphism of the
      star product $\stars$ of standard ordered type which again is a
      formal series in $\lambda$ whose coefficients consist of differential
      operators. Hence it suffices to prove the desired
      identity $\Cs^{-1} N^2 \mathcal C = \id$ locally on functions
      in $C_{pp}^\infty(T^*Q)$. Take again the above neighbourhood $U$.
      By proposition \ref{AnalProp} every such function can be written
      as a finite linear combination over terms consisting of star products
      of functions belonging to
      $C^\infty_{pp,0}(T^*Q)[\lambda] + C^\infty_{pp,1}(T^*Q)[\lambda]$.
      By the automorphism property it thus suffices to check our
      identity on these functions at most linear in the momenta.
      An easy computation using lemma \ref{CSantiaut} shows that this
      is the case.

\item This is straight forward: using the fact that $N^2 \mathcal C$
      is an involutive anti-linear anti-automorphism of $\stars$ it
      follows that complex conjugation is an involutive anti-linear
      anti-automorphism of $\starw$. Moreover since all constructions
      depend only on the combination $\frac{\lambda}{\im}$ it follows
      that the parity transformation generated by
      $\lambda \mapsto - \lambda$ is an anti-automorphism which shows
      that $\starw$ is indeed of Weyl type.

\item This is a straight forward computation.
\end {enumerate}
\end {PROOF}

{\bf Remark:} The two star products $\starw$ and $\starf$ are not the same 
in general as we shall see evaluating both star 
products on functions linear in the momenta: Let $X, Y \in \Gamma (TQ)$ 
and consider $\widehat X, \widehat Y \in C^\infty_{pp,1} (T^*Q)$. 
Then one gets by direct calculation 
$N \widehat X = \widehat X + \frac{\lambda}{2\im} \Div_\alpha X$
as well as 
\BEQ {DeltahatXY}
    \Delta (\widehat X \widehat Y)
    =
    \widehat X \Div_\alpha Y + \widehat Y \Div_\alpha X + 
    \widehat {{\nabla_0}_X Y} + \widehat {{\nabla_0}_Y X}
\EEQ    
\BEQ {DeltaDeltaYX}
    \Delta^2 (\widehat X \widehat Y) =
    2 (\Div_\alpha X) (\Div_\alpha Y)
    + \Lie_X (\Div_\alpha Y) + \Lie_Y (\Div_\alpha X)
    + \Div_\alpha \left({\nabla_0}_X Y \right) 
    + \Div_\alpha \left({\nabla_0}_Y X \right) .
\EEQ    
Using theorem \ref {StandardRepTheo} one obtains 
$\srep (N\widehat X) = \frac{\lambda}{\im} (\Lie_X + \frac{1}{2} 
\Div_\alpha X)$ and 
\BEQ {SrepYX}
    \srep \left(\widehat X \widehat Y\right) 
    = \left(\frac{\lambda}{\im}\right)^2 \Lie_X \Lie_Y 
    - \frac{\lambda}{\im} \srep \left(\widehat{{\nabla_0}_X Y}\right) .
\EEQ
This enables us to calculate $\srep (N\widehat X) \srep (N\widehat Y)$
and using the injectivity of $\srep$ on functions polynomial in the 
momenta (corollary \ref {SRepInjCor}) and the representation property
of $\srep$ we obtain 
\begin {eqnarray*}
    \lefteqn {(N\widehat X) \stars (N\widehat Y)} \\
    & = & \widehat X \widehat Y 
    + \frac{\lambda}{\im} \left(\widehat{{\nabla_0}_X Y} 
    + \frac{1}{2} \widehat Y \Div_\alpha X 
    + \frac{1}{2} \widehat X \Div_\alpha Y \right)
    + \frac{1}{2} \left(\frac{\lambda}{\im}\right)^2
    \left( \Lie_X (\Div_\alpha Y) 
    + \frac{1}{2} (\Div_\alpha X)(\Div_\alpha Y) \right) .
\end {eqnarray*}
Finally we calculate $\widehat X \starw \widehat Y = N^{-1} ((N\widehat X) 
\stars (N\widehat Y))$ and obtain
\BEQ {XstarwY}
    \widehat X \starw \widehat Y =
    \widehat X \widehat Y 
    + \frac{\im\lambda}{2} \left\{ \widehat X, \widehat Y \right\}
    + \left(\frac{\im\lambda}{2}\right)^2 
    M^\TinyW_2 (\widehat X, \widehat Y)
\EEQ
where 
\BEQ {MZweiW}
    M^\TinyW_2 (\widehat X, \widehat Y) =
    \frac{1}{2} \Big(
    \Lie_X (\Div_\alpha Y) + \Lie_Y (\Div_\alpha X) 
    - \Div_\alpha \left({\nabla_0}_X Y \right)   
    - \Div_\alpha \left({\nabla_0}_Y X \right) \Big)
\EEQ
observing that 
$\Lie_X (\Div_\alpha Y) - \Lie_Y (\Div_\alpha X) 
- \Div_\alpha ({\nabla_0}_X Y) + \Div_\alpha ({\nabla_0}_Y X) = 0$ 
due to $d\alpha (X, Y) = - \tr R_0 (X, Y)$. Writing the coordinate
expression $X^k_{|l} Y^l_{|k}$
for ${\rm trace}(Z\mapsto {\nabla_0}_{{\nabla_0}_ZX}Y)$ to avoid
clumsy notation the operator $M^\TinyW_2$ can be simplified to
\BEQ {MZweiWLocal}
    M^\TinyW_2 (\widehat X, \widehat Y) =
    -  X^k_{|l} Y^l_{|k}
    - \frac{1}{2} \Big(\Ric_0 (X, Y) + \Ric_0 (Y, X)
    - \left({\nabla_0}_X \alpha\right) (Y) 
    - \left({\nabla_0}_Y \alpha\right) (X) \Big)
\EEQ
where ${\,}_{|k}$ denotes the covariant derivative with respect to 
$\partial_{q^k}$ and $\Ric_0$ denotes the Ricci tensor of $\nabla_0$. 
On the other hand we compute the Fedosov star 
product $\starf$ of $\widehat X$ and $\widehat Y$ using 
\cite [Theorem 3.4] {BW96a} and obtain here the following expression
for the second order term $M_2^\TinyF (\widehat X, \widehat Y)$:
\BEQ {MZweiFLocal}
    M^\TinyF_2 (\widehat X, \widehat Y) 
    =
    - X^k_{|l} Y^l_{|k}~~.
\EEQ
For general manifolds $Q$ with torsion-free connection $\nabla_0$ these
two star
products do {\ not} coincide for any choice of $\alpha$: for example,
let
$Q$ be equal to $S^2$ with the standard metric. Its Levi-Civita connection
$\nabla_0$ is clearly unimodular whence every possible $\alpha$ is a
closed one-form, hence exact ($\alpha=d\phi$) since the the first de Rham
cohomology group of the two-sphere vanishes. If the two expressions
(\ref{MZweiWLocal}) and (\ref{MZweiFLocal}) were the same for any $X,Y$
then the Ricci tensor would be equal to the second covariant derivative
of $\phi$. In particular upon contracting with the inverse metric we
would obtain that the Laplacian of $\phi$ were equal to a positive
multiple of the scalar curvature of $S^2$ which is positive and
constant. But this differential equation has no smooth solution since
the integral (with respect to the Riemannian volume) over $S^2$ of the
Laplacian of $\phi$ would vanish as opposed to the integral over $S^2$
of a positive constant.

\section {Example: The cotangent bundle of a Lie group}
\label {LieSec}

This section shall be dedicated to finding explicit formulae for the
homogeneous Fedosov star product $\stars$ of standard ordered type
on the cotangent bundle of a connected Lie group $G$, that is equipped
with the natural torsion-free connection defined by
${\nabla_0}_U V :=\frac{1}{2}[U,V]$ for left-invariant vector fields
$U,V$ on $G$, obtained using the standard representation
$\srep$ for functions polynomial in the momentum variables.
To express the bidifferential operators
defining the star product we make use of the natural set of
basis sections
in the tangent bundle of $T^*G$ given by $Y_i$ resp. $Z^j$ that are the
horizontal resp. vertical lifts of a basis of left-invariant vector
fields $X_i$ resp. the dual left-invariant one-forms $\theta^j$ with
respect to the flat connection on $G$ (see definition \ref {LiftDef})
which is defined by $\widetilde\nabla_U V =0$ for left-invariant vector
fields $U$, $V$. Instead of using Darboux coordinates it will be
convenient to use natural fibre-variables $P_i~(1\leq i \leq \dim(G))$
given by
$P_i : T^*G \rightarrow {\mathbb R} : \alpha_g \rightarrow \alpha_g(X_i)$
for $\alpha_g \in T_g^*G$.
The obtained expression for $\stars$ shall moreover be related to a
star product $\starg$ of Weyl type by means of an
equivalence transformation $N$ and we shall prove that
$\starg$ coincides (up to a rescaling of the formal
parameter) with the star product obtained by Gutt in \cite{Gut83}
using cohomological methods instead of our purely algebraic approach.
As a first step we shall find an expression for the star product $\stars$
of polynomial functions in the momenta on $T^*G$ that are invariant under
the canonical lift $T^*(l_g)$ of the left-translations
$l_g : G \rightarrow G$ to $T^*G$. Since those polynomials are
generated by functions of the form
$\widehat{U^{\vee k}}$ for left-invariant vector fields $U$ on $G$,
we may restrict our calculations to functions
\[
    e_U : \alpha_g \mapsto
    \sum\limits_{r=0}^{\infty}
    \frac{1}{r!}\left(U(\alpha_g)\right)^r
    \quad
    \in C^\omega_p (T^*G)[[\lambda]].
\]
\begin {LEMMA} \label{g1}
Let $U,V$ be left-invariant vector fields on $G$ then we have
the formula \cite {Dri83}, \cite{Gut83}
\BEQ {DrinfeldStar}
    e_U \stars e_V =
    e_{\frac{\im}{\lambda}
       H \big(\frac{\lambda}{\im} U, \frac{\lambda}{\im} V\big)},
\EEQ
at which $H$ denotes the Baker-Campbell-Hausdorff series. Moreover this
equation uniquely determines bidifferential operators $M_r^{inv}$
defined by
\BEQ {MinvDef}
    e_U \stars e_V =:
    \sum\limits_{r=0}^\infty \left(\frac{\lambda}{\im}\right)^r
    M_r^{inv}(e_U,e_V),
\EEQ
that are of order $r$ in every argument, homogeneous of degree $-r$
and only containing derivatives with respect to vertical
directions and multiplications with polynomials in the momenta
on $T^*G$, that are
invariant under $T^*(l_g)$.
\end {LEMMA}
\begin {PROOF}
We just have to notice that $\srep (e_U) \chi$ is given by
$\exp(\circ \frac{\lambda}{\im}U) \chi$ for $\chi \in C^\infty(G)$,
at which $\circ$ denotes the composition of the left-invariant vector
fields viewed as differential operators on $C^\infty(G)$.
Using this fact yields
\begin {eqnarray*}
    \srep (e_U \stars e_V) \chi
    & = & \left(\srep (e_U)\circ \srep(e_V)\right) \chi =
          \left(\exp\left(\circ \frac{\lambda}{\im} U \right)
          \circ \exp\left(\circ \frac{\lambda}{\im} V \right)\right) \chi \\
    & = & \exp\left(\circ \frac{\lambda}{\im}\frac{\im}{\lambda}
          H\left(\frac{\lambda}{\im}U,
          \frac{\lambda}{\im}V\right)\right) \chi
          = \srep\left(
          e_{\frac{\im}{\lambda}H
          \big(\frac{\lambda}{\im}U, \frac{\lambda}{\im}V \big)}\right)
          \chi.
\end {eqnarray*}
Now both arguments $e_U \stars e_V$ and
$e_{\frac{\im}{\lambda} H\big(\frac{\lambda}{\im} U, \frac{\lambda}{\im}
V\big)}$ are elements in the submodule
$C^\omega_p (T^*Q) [[\lambda]]$ and hence corollary \ref {SRepInjCor}
ensures that they coincide.
The assertions about the bidifferential operators $M_r^{inv}$ are obvious
consequences of properties of the Baker-Campbell-Hausdorff series resp.
of proposition \ref {DstdProp}.
\end {PROOF}
\begin {PROPOSITION} \label {g2}
The star product $\stars$ of two functions
$f, g \in C^\infty(T^*G)$ is given by
\BEQ {LieStdStar}
    f \stars g =
    \sum\limits_{r=0}^\infty
    \left(\frac{\lambda}{\im}\right)^r
    \sum\limits_{t=0}^r \frac{1}{t!}
    M_{r-t}^{inv}(Z^{j_1} \ldots Z^{j_t}f, Y_{j_1} \ldots Y_{j_t} g).
\EEQ
By lemma \ref {g1} it is obvious that this star product is of Vey type.
\end {PROPOSITION}
\begin {PROOF}
First we notice that it is sufficient to prove (\ref {LieStdStar})
for functions polynomial in the momenta. Using
$\srep (\widehat{S}) \chi = \frac{1}{l!} \left(\frac{\lambda}{\im}\right)^l
S^{j_1 \ldots j_l} X_{j_1} \ldots X_{j_l} \chi$ for
$\widehat{S} \in C^\infty_{pp,l}(T^*G)$ and an analogous formula for
$\widehat{T} \in C^\infty_{pp,k}(T^*G)$,which are obtained from
$D^k_0 \chi = (\theta^l \vee {\nabla_0}_{X_l})^k \chi =
X_{i_1} \ldots X_{i_k} \chi ~\theta^{i_1} \vee \ldots \vee \theta^{i_k}$
and the symmetry of $S$ resp. $T$,
we get for $\chi \in C^\infty (G)$
\begin{eqnarray*}
    \srep(\widehat{S} \stars \widehat{T}) \chi
    & \stackrel{\textrm{(a)}}{=} &
          \sum\limits_{t=0}^l
          \frac{1}{k!t!(l-t)!}\left(\frac{\lambda}{\im}\right)^t
          S^{j_1 \ldots j_l}
          (X_{j_1} \ldots X_{j_t}T^{i_1 \ldots i_k})
          \srep (P_{j_{t+1}} \cdots P_{j_l}
          \stars P_{i_1} \cdots P_{i_k}) \chi \\
    & \stackrel{\textrm{(b)}}{=} &
          \srep \left(\sum\limits_{r=0}^\infty
          \left(\frac{\lambda}{\im}\right)^r \sum\limits_{t=0}^l
          \frac{1}{t!(l-t)!}
          \left(\frac{\lambda}{\im}\right)^t
          M_r^{inv}\left(
          \pi^*S^{j_1 \ldots j_l}P_{j_{t+1}} \cdots P_{j_l},
          Y_{j_1} \ldots Y_{j_t}\widehat{T}\right)\right) \chi \\
    & \stackrel{\textrm{(c)}}{=} &
          \srep\left(\sum\limits_{r=0}^\infty \left(
          \frac{\lambda}{\im}\right)^r \sum\limits_{t=0}^r
          \frac{1}{t!} M_{r-t}^{inv}
          \left(Z^{j_1} \ldots Z^{j_t}
          \widehat{S},Y_{j_1} \ldots Y_{j_t} \widehat{T}\right)\right)
          \chi.
\end{eqnarray*}
Equation (a) is a consequence of the Leibniz rule and
$\srep(P_{i_1} \cdots P_{i_k})\chi = \left(\frac{\lambda}{\im} \right)^k
\frac{1}{k!} \sum\limits_{\sigma \in S_k} X_{i_{\sigma (1)}} \cdots
X_{i_{\sigma (k)}} \chi$.
In (b) we used lemma \ref{g1} and the fact that
$\srep(\pi^*\psi F)= \psi \srep(F)$
for $\psi \in C^\infty(G)$ and the $\pi$-relatedness of the
vector fields $X_i$ and
$Y_i$. In (c) the sum over $t$ was extended to $\infty$ since the terms
$Z^{j_1} \ldots Z^{j_t}\widehat{S}$ vanish for $t>l$ due to lemma
\ref {WichtigLem}. Then (\ref {LieStdStar}) follows by the injectivity
of $\srep$ restricted to $C^\infty_{pp} (T^*G)$ proved in
corollary \ref {SRepInjCor}.
\end {PROOF}

Now according to the last section we consider a Weyl ordered star
product defined by
\BEQ {GuttStarDef}
    f \starg g :=
    N^{-1}\left(N f \stars N g \right)
    \textrm{ and }
    N:=\exp \left( \frac{\lambda}{2\im} \Delta \right)
\EEQ
with $\Delta = Y_i Z^i + \frac{1}{2} C^l_{li} Z^i + \pi^*(\alpha_i)
Z^i$ (which is easily computed using equation (\ref{CovDeltaDef}))
denoting by $\alpha_i$ the components of a one-form fulfilling
$d\alpha = -\tr R_0=0$ due to the particular choice of the
connection $\nabla_0$ where 
$C^k_{ij} := \theta^k ([X_i, X_j])$ are the structure constants of the 
Lie algebra $\mathfrak g$ of $G$. Choosing 
$\alpha = (t-\frac{1}{2}) C^l_{li} \theta^i$ for $t \in \mathbb R$ 
it is easily verified that
$d\alpha=0$ and $N = N_0 A_t$ with $N_0:= \exp \left(
\frac{\lambda}{2\im} Y_i Z^i\right)$ and $A_t:= \exp
\left(\frac{\lambda}{2\im}t d_\TinyS \right)$ with $d_\TinyS:=
C^l_{li} Z^i$, since $Y_i Z^i$ and $d_\TinyS$ commute. For any
choice of $t \in \mathbb R$ we shall prove the coincidence of
$\starg$ with the one constructed in \cite{Gut83}.
\begin {LEMMA} \label{LeftinvPolyCoinc}
For $\alpha=(t-\frac{1}{2}) C^l_{li} \theta^i ~(t \in \mathbb R)$
and $N$ defined as above we have
\BEQ {coinceq}
    e_U \starg e_V  = e_U \stars e_V
    \textrm{ for all left-invariant vector fields $U,V$ on G.}
\EEQ
Moreover $A_t$ is a one-parameter-group of automorphisms of $\stars$,
i.~e. $A_t A_s = A_{t+s}$ and $A_0 = \id$ and
\BEQ {Automorph}
    A_t\left( f \stars g \right)= (A_t f) \stars (A_t g) \quad
    \forall f, g \in C^\infty (T^*G)[[\lambda]], \forall t \in \mathbb R,
\EEQ
implying that $d_\TinyS$ is a derivation of $\stars$.
Moreover this implies that all choices of $t \in \mathbb R$ lead to the
same star product $\starg$.
\end {LEMMA}
\begin {PROOF}
Since $N=N_0 A_t$ we prove the first assertion in two steps using
$N_0$ and $A_t$ separately.
Obviously we have $N_0 e_U = e_U$ and therefore using equation
(\ref{DrinfeldStar})
\[
    {N_0}^{-1}\left( (N_0 e_U) \stars (N_0 e_V) \right) =
    {N_0}^{-1}\left(e_U \stars e_V\right) =
    {N_0}^{-1} e_{\frac{\im}{\lambda}
    H \big(\frac{\lambda}{\im} U, \frac{\lambda}{\im} V\big)}=
    e_{\frac{\im}{\lambda}
    H \big(\frac{\lambda}{\im} U, \frac{\lambda}{\im} V\big)}=
    e_U \stars e_V.
\]
For the second part we compute $\srep(A_t e_U) =
\exp\left(\frac{\lambda}{2\im}t C^l_{li} U^i\right) \srep (e_U)$
yielding
\[
    \srep\left( (A_t e_U) \stars (A_t e_V) \right) =
    \exp\left(\frac{\lambda}{2\im} t C^l_{li} (U^i + V^i)\right)
    \srep (e_U \stars e_V).
\]
On the other hand we have again by equation (\ref{DrinfeldStar})
\[
    \srep\left( A_t (e_U \stars e_V) \right) =
    \exp\left(\frac{\lambda}{2\im} t C^l_{li} (U^i + V^i)\right)
    \srep \left(e_{\frac{\im}{\lambda}
    H \big(\frac{\lambda}{\im} U, \frac{\lambda}{\im} V\big)}\right)=
    \exp\left(\frac{\lambda}{2\im} t C^l_{li} (U^i + V^i)\right)
    \srep (e_U \stars e_V)
\]
using the shape of the Baker-Campbell-Hausdorff series and the fact that
$C^l_{li} [W,X]^i =0$ for any left-invariant vector fields $W,X$ on $G$.
By corollary \ref{SRepInjCor} this implies
$A_t^{-1}\left( (A_t e_U) \stars (A_t e_V) \right) = e_U \stars e_V$.
Since $d_\TinyS$ only contains derivatives with respect to vertical
directions this implies that $A_t$ is even an automorphism of $\stars$ for
all $t \in \mathbb R$, proving the second assertion. The derivation property
of $d_\TinyS$ is an immediate consequence of (\ref {Automorph}).
\end {PROOF}

Observe that an equation analogous to (\ref{Automorph}) does not hold
for $N_0$, instead we have the following lemma:
\begin {LEMMA} \label{g4}
For any $f \in C^\infty(T^*G)$ and $\chi \in C^\infty(G)$ we have
\BEQ {pichiGuttF}
    \pi^*\chi \starg f = \sum\limits_{r=0}^\infty \frac{1}{r!}
    \left(\frac{\im \lambda}{2} \right)^r
    \pi^*(X_{i_1} \ldots X_{i_r} \chi) Z^{i_1} \ldots Z^{i_r} f .
\EEQ
\end {LEMMA}
\begin {PROOF}
The proof is a lengthy but straight forward computation using
proposition \ref{g2}
and the fact that all terms containing $Z$-derivations applied
to $\pi^*\chi$ vanish due to lemma \ref {WichtigLem}.
\end {PROOF}
\begin {LEMMA} \label{g5}
For any left-invariant vector fields $U,V$ on $G$ we have
\BEQ {UGutteV}
    \widehat{U} \starg e_V =
    e_V \widehat{\frac{[\frac{\lambda}{\im}V,\cdot ]}
    {\exp([\frac{\lambda}{\im}V,\cdot ])-1}U}.
\EEQ
Expressing the term $W_r$ of order $\left(\frac{\lambda}{\im}\right)^r$
in components this reads
\BEQ {Wr}
    W_r(\widehat{U},e_V) =
    \frac{(-1)^r}{r!}B_r P_l (Z^j \widehat{U})C^{j_1}_{jk_1} \cdots
    C^{j_{r-1}}_{j_{r-2}k_{r-1}}C^l_{j_{r-1}k_r}(Z^{k_1} \ldots Z^{k_r} e_V)
\EEQ
where $B_r$ denotes the $r$th Bernoulli number.
\end {LEMMA}
\begin {PROOF}
By lemma \ref{LeftinvPolyCoinc}
$e_U \stars e_V = e_U \starg e_V$ and therefore we have
\[
    \widehat{U} \starg e_V
    = \left.\frac{d}{dt}\right|_{t=0} e_{tU} \starg e_V
    = \left.\frac{d}{dt}\right|_{t=0}
      e_{\frac{\im}{\lambda}H(t \frac{\lambda}{\im}U,
      \frac{\lambda}{\im}V)}
    = e_V \widehat{\frac{[\frac{\lambda}{\im}V,\cdot ]}
      {\exp([\frac{\lambda}{\im}V,\cdot ])-1}U}
\]
using that for $a,b \in \mathfrak{g}$ the following is valid
\[
    \left.\frac{d}{dt} \right|_{t=0} H(ta,b) =
    \frac{\ad b}{\exp(\ad b)-1}a.
\]
Then equation (\ref {Wr}) is a direct consequence of the Taylor series
expansion of $\frac{x}{e^x -1}$.
\end {PROOF}

Collecting our results we get the following proposition.
\begin {PROPOSITION} \label {g6}
The star product $\starg$ being related to the Fedosov star product 
of standard ordered type $\stars$ by means of the equivalence 
transformation $N$ coincides with the product constructed in \cite{Gut83}.
\end {PROPOSITION}
\begin {PROOF}
The assertion follows from lemma \ref{g4} and lemma \ref{g5} and the
theorem proven in chapter 4 in \cite{Gut83}.
\end {PROOF}


\appendix

\section {Homogeneous connections on cotangent bundles}
\label {ConnApp}

We now shall briefly recall some well-known basic definitions
concerning horizontal and vertical lifts of vector fields,
homogeneous connections and normal Darboux
coordinates on cotangent bundles.
\begin {DEFINITION} \label {LiftDef}
Let $\nabla_0$ be a connection on $Q$ then we consider the
connection mapping $K : T(T^*Q) \rightarrow T^*Q$
of $\nabla_0$ defined by
\[
    K(\dot{\alpha}(0)) :={\nabla_0}_
    {\partial_t} \alpha \mid_{t=0}
\]
for a curve $\alpha$ in $T^*Q$. It turns out that
$(T \pi \times K) : T(T^*Q) \rightarrow TQ \oplus T^*Q$
is a fibrewise isomorphism. Because of this fact the notion of
horizontal resp. vertical lifts with respect to
$\nabla_0$ is well-defined by: $X^h \in \Gamma(T(T^*Q))$
is called the horizontal lift of $X \in \Gamma(TQ)$ iff
\BEQ {HLiftDef}
    T\pi X^h =X \circ \pi \textrm{ and } K(X^h)=0
\EEQ
resp. $\beta^v \in \Gamma(T(T^*Q))$ is called the vertical lift of
$\beta \in\Gamma(T^*Q)$ iff
\BEQ {VLiftDef}
    K(\beta^v)=\beta \circ \pi \textrm{ and } T\pi \beta^v =0.
\EEQ
\end {DEFINITION}
Working in a local bundle Darboux chart
$(q^1, \ldots, q^n, p_1, \ldots, p_n)$ one finds that
\BEQ {LiftVectDef}
    X^h = \left(\pi^*X^i\right) \partial_{q^i} + \pi^*(X^k
          \Gamma^j_{ki}) p_j \partial_{p_i}
    \qquad
    \mbox{ and }
    \qquad
    \beta^v = \left(\pi^* \beta_i\right)
    \partial_{p_i},
\EEQ
at which $\Gamma^j_{ki}$ denotes the Christoffel symbols
of $\nabla_0$ and $X^i$ resp. $\beta_i$ the
components of $X$ resp. $\beta$.
The following lemma should be well-known and is crucial throughout
this paper. It can be proved easily in a local bundle chart:
\begin {LEMMA} \label {WichtigLem}
Let $T \in \Gamma (\bigotimes^l T^*(T^*Q))$  where $l \ge 0$
be a homogeneous tensor field of degree $k$, i.~e. $\Lie_\xi T = kT$.
Then $k\in \mathbb Z$ and $T=0$ if $k<0$ and $T = \pi^*\tilde T$ with
$\tilde T \in \Gamma (\bigotimes^l T^*Q)$ iff $k = 0$.
Moreover $[\Lie_\xi , \Lie_{\beta^v}] = - \Lie_{\beta^v}$ and
$[\Lie_\xi, i_m (\beta^v)]T = - i_m (\beta^v)T$ where $i_m$ denotes the
substitution into the $m$th argument of $T$. Hence $\Lie_\xi T = kT$ 
implies $\Lie_\xi (\Lie_{\beta^v} T) = (k-1) T$ as well as
$\Lie_\xi (i_m (\beta^v) T) = (k-1) i_m (\beta^v) T$.
\end {LEMMA}
\begin {DEFINITION} \label {HomConnDef}
A connection $\nabla$ on $T^*Q$ is said to be homogeneous iff
$\Lie_\xi \nabla = 0$, i.~e.
$\Lie_\xi \nabla_X Y - \nabla_{\Lie_\xi X} Y - \nabla_X \Lie_\xi Y = 0$
for all $X, Y \in \Gamma (T(T^*Q))$.
\end {DEFINITION}
\begin {DEFINITION} \label {LiftedConnDef}
The remark about $(T\pi \times K)$ being a fibrewise isomorphism justifies
defining a connection on $T^*Q$ by disposing of
$\nabla^0$ for $X,Y \in \Gamma(TQ), \beta,\gamma \in \Gamma(T^*Q)$
in the following manner:
\begin{eqnarray}
    \left.\nabla^0_{X^h}Y^h \right|_{\alpha_q}
    & := & \left.\left({\nabla_0}_X Y\right)^h \right|_{\alpha_q}
           + \alpha_q \left.\left(
           \frac{1}{2}R(X, Y) \underline{\,\cdot\,} \,
           - \frac{1}{6}
           \left(R (\underline{\,\cdot\,}, X) Y +
                 R (\underline{\,\cdot\,}, Y) X\right)
           \right)^v \right|_{\alpha_q} \\
    \left.\nabla^0_{X^h} \beta^v \right|_{\alpha_q}
    & := & \left.\left({\nabla_0}_X \beta\right)^v \right|_{\alpha_q} \\
    \left.\nabla^0_{\beta^v}X^h \right|_{\alpha_q}
    & := & \left.\nabla^0_{\beta^v} \gamma^v \right|_{\alpha_q}
      := 0.
\end{eqnarray}
\end {DEFINITION}
It is straight forward to check that $\nabla^0$ is homogeneous,
torsion-free and symplectic. A simple calculation using a bundle Darboux
chart yields the following local expressions for the Christoffel symbols
${\Gamma^0}^{x^k}_{x^i x^j}$ and the components of the
curvature tensor
${R^0}^{x^k}_{x^i x^j x^l} =
dx^k (\nabla^0_{\partial_{x^j}} \nabla^0_{\partial_{x^l}} \partial_{x^i}
-\nabla^0_{\partial_{x^l}} \nabla^0_{\partial_{x^j}} \partial_{x^i})$
of the connection $\nabla^0$:
\[
    {\Gamma^0}^{q^k}_{q^i q^j} =
    - {\Gamma^0}^{p_j}_{q^i p_k} =
    - {\Gamma^0}^{p_i}_{p_k q^j}
    = \pi^* \Gamma^k_{ij}
\]
\[
    {\Gamma^0}^{p_k}_{q^i q^j}
    = \frac{p_a}{3} \pi^* \left(
    2 \Gamma^a_{js} \Gamma^s_{ki} -
    \partial_{q^j} \Gamma^a_{ki} +
    \textrm{cycl.}(ijk)
    \right)
\]
\[
    {R^0}^{q^l}_{q^k q^i q^j} =
    -{R^0}^{p_k}_{p_l q^i q^j} = \pi^* R^l_{kij}
    \qquad
    {R^0}^{p_l}_{q^k q^i p_j} =
    \frac{1}{3}\pi^*\left(R^j_{lki} + R^j_{kli}\right)
\]
\[
    {R^0}^{p_i}_{q^j q^k q^l}
    = \frac{p_a}{3} \pi^* \left(
      R^a_{jlk|i}
      - 3 \Gamma^a_{is} R^s_{jlk} - \Gamma^a_{ls} R^s_{ijk} +
      \Gamma^a_{ks} R^s_{ijl} + (i \leftrightarrow j)
   \right)
\]
All other combinations vanish and $\ldots_{|i}$ denotes the covariant
derivative with respect to $\partial_{q^i}$.
At this instance we want to refer to a question that was brought to
our interest by M.~Cahen namely whether the Ricci tensor $\Ric^0$
corresponding to $R^0$ defined by
$\Ric^0(X,Y):=\tr (Z \rightarrow R^0(Z,X)Y)$ for
$X,Y \in \Gamma(T(T^*Q))$ enjoys the property that
$(\nabla^0_X \Ric^0)(Y,Z) + cycl.(X,Y,Z) =0$ for all
$X,Y,Z \in \Gamma(T(T^*Q))$. It turns out by a direct but lengthy
calculation that this is the fact iff the Ricci tensor $\Ric_0$
corresponding to $R_0$ does so with respect to $\nabla_0$.
In the next definition we shall briefly
explain the concept of a connection on $Q$ that is induced by a
connection on $T^*Q$.
\begin {DEFINITION} \label {InducedConnDef}
Let $\overline{\nabla}$ be a torsion-free homogeneous
connection on $T^*Q$. Choose any connection $\nabla^{\mbox{\rm\tiny Q}}$
on $Q$ and define for $X,Y \in \Gamma(TQ)$:
\[
    Ti \left({\nabla_0}_X Y\right) :=
    \overline{\nabla}_{X^h} Y^h \circ i
\]
where the horizontal lift is taken with respect to
$\nabla^{\mbox{\rm\tiny Q}}$. Then ${\nabla_0}_X Y$ is a
well-defined vector field on $Q$, and does not depend on the choice of
$\nabla^{\mbox{\rm\tiny Q}}$. We refer to $\nabla_0$
as the connection induced by $\overline{\nabla}$.
\end {DEFINITION}
Using this notion we can give a characterization of
$\nabla^0$ as follows:
\begin {PROPOSITION} \label {iaXDeltaInvRProp}
Let $\nabla_0$ be a torsion-free connection on
$Q$. Among all homogeneous, torsion-free, symplectic connections
$\overline{\nabla}$ on $T^*Q$ which induce $\nabla_0$ on $Q$
the connection $\nabla^0$ is uniquely characterized
by the condition
\BEQ {iaXDeltaInv}
    i_a(X) \delta^{-1} \overline{R} = 0
\EEQ
for all vertical vector fields $X$ on $T^*Q$ where
$\overline R := \frac{1}{4} \omega_{it} \overline R^{t}_{jkl}
dx^i \vee dx^j \otimes dx^k \wedge dx^l$.
\end {PROPOSITION}
\begin {PROOF}
The proof mainly relies on the fact that there is a uniquely determined
tensor field $B \in \Gamma(TQ \otimes \bigvee^3 T^*Q)$ such that
\[
    \overline{\nabla}_{X^h} Y^h \mid_{\alpha_q} =
    {\nabla^0}_{X^h} Y^h \mid_{\alpha_q}
    +\alpha_q \left(B(X,Y,\underline{\,\cdot\,})\right)^v \mid_{\alpha_q}
\]
for all $X,Y \in \Gamma(TQ)$ the horizontal lifts being taken with
respect to $\nabla_0$.
Then the assertion follows by comparing $\overline{R}$ to
$R^{\mbox{\rm\tiny 0}}$ and the fact that (\ref{iaXDeltaInv}) is
fulfilled iff $B=0$.
\end {PROOF}

Finally we shall prove the following three lemmata:
\begin {LEMMA} \label {NormalDarbouxLem}
Let $(q^1, \ldots, q^n)$ be a normal chart around
$q_0$ with domain $D^q$ of $Q$ with respect to $\nabla_0$ then for any
$\alpha_q \in \pi^{-1}(D^q)$ we have
\[
    \left.{\Gamma^0}_{ijk}(x) x^i x^j x^k
    \right|_{\alpha_q}
    :=
    \left.\omega_{il}
    {\Gamma^0}^l_{jk}(x) x^i x^j x^k
    \right|_{\alpha_q}=0,
\]
where $(x^1,\ldots,x^{2n}) = (q^1,\ldots,q^n,p_1,\ldots,p_n)$. Hence
the bundle Darboux coordinates are normal Darboux coordinates with
respect to $\nabla^0$ around $\alpha_{q_0}$.
\end {LEMMA}
\begin {PROOF}
By direct calculation the assertion follows easily using
the fact that $(q^1, \ldots, q^n)$ are normal coordinates on $Q$
and the local expressions for ${\Gamma^0}^{x^k}_{x^i x^j}$ as stated
above.
\end {PROOF}
\begin {LEMMA} \label {TrRExactLem}
Let $M$ be a differentiable manifold and $\nabla$ a connection on $M$
then the trace of the curvature tensor $R$ is exact i.~e.
\[
    (\tr R)(X,Y) :=dx^i \left(R(X,Y) \partial_i \right) =
    d \alpha (X,Y) \quad \forall X,Y \in \Gamma(TM)
\]
for an $\alpha \in \Gamma(T^*M)$. If $\tr R=0$ then the connection
$\nabla$ is called unimodular.
\end {LEMMA}
\begin {PROOF}
First observe that the Levi-Civita connection
$\nabla^{\mbox{\rm\tiny LC}}$ of a Riemannian metric
is unimodular. Now since every manifold admits such a connection we
shall compare this one to $\nabla$.
Let $S \in \Gamma(TM \otimes T^*M \otimes T^*M)$
be the uniquely determined tensor field such that
${\nabla^{\mbox{\rm\tiny LC}}}_X Y =\nabla_X Y + S(X,Y)$.
Straight forward computation using this equation and
$\tr R^{\mbox{\rm\tiny LC}} = 0$ yields
\[
    (\tr R)(X,Y) =
    - \tr \left(Z \mapsto \left(
    S(T(X, Y), Z) + S(X, S(Y, Z)) - S(Y,S(X, Z))
    + (\nabla_X S)(Y,Z) - (\nabla_Y S)(X,Z) \right)\right)
\]
where $T$ is the torsion of $\nabla$.
Since $\tr \left(Z \rightarrow S(X,S(Y,Z)) - S(Y,S(X,Z))\right) = 0$
and since $\nabla$ commutes with contractions we get
$(\textsf{tr}R)(X,Y) = \left((\nabla_X \alpha)(Y) -
(\nabla_Y \alpha)(X)\right) + \alpha (T(X, Y)) = d\alpha (X,Y)$
for $\alpha (Y) := - \tr \left( Z \rightarrow S(Y,Z)\right)$.
\end {PROOF}

The last lemma should be well-known and is used in
theorem \ref {ClassFTTheo}:
\begin{LEMMA} \label{TaylorLem}
Let $\nabla$ be a torsion-free connection on a connected manifold $M$.
Considering normal coordinates $(x^1,\ldots ,x^n) (\dim(M) = n)$
around an arbitrary point $q \in M$
we have the following identity for the Taylor series of a function
$f \in C^\infty(M)$ with respect to $\nabla$
\[
    \left.\tau_0 (f)\right|_q :=\sum\limits_{r=0}^{\infty}
    \left. \frac{1}{r!} D^r f \right|_q
    =\sum\limits_{r=0}^{\infty} \frac{1}{r!} \left.\frac
    {\partial^r f}{\partial x^{i_1} \cdots \partial x^{i_r}}
    \right|_q
    dx^{i_1} \vee \ldots \vee dx^{i_r}
\]
where $D := dx^k \vee \nabla_{\partial_{x^k}}$.
\end{LEMMA}
\begin{PROOF}
For the proof one has to observe that in any local coordinates
$D^r S = (dx^i \vee \Lie_{\partial_{x^i}} - dx^i \vee dx^k
\vee \Gamma_{ik}^j i_s(\partial_{x^j}))^r S $ where $S$ is
an arbitrary symmetric covariant tensor field on $M$. By induction
this is equal to $(dx^i \vee \Lie_{\partial_{x^i}})^r S
+ \sum_{a=1}^{N_r} C_a (K_a (S))$ where $N_r$ is some integer,
$K_a$ are certain differential operators and
$C_a = dx^{i_1} \vee \cdots \vee dx^{i_{s+2}}
(\partial_{x^{i_1}} \cdots \partial_{x^{i_s}} \Gamma^j_{i_{s+1} i_{s+2}})
D_j$ where either
$D_j = \Lie_{\partial_{x^j}}$ or $D_j = i_s (\partial_{x^j})$ and
$s \ge 0$. Finally note that in normal coordinates around $q$
we have
$0=(\partial_{x^{i_1}} \cdots \partial_{x^{i_r}} \Gamma_{kl}^a)(q) dx^k
\vee dx^l \vee dx^{i_1} \vee \ldots \vee dx^{i_r}$ at $q$
which is obtained by
$r$-fold differentiation (with respect to the affine parameter)
of the geodesic equation for a geodesic emanating at $q$.
Hence $C_a$ equals $0$ at $q$ which proves the lemma after having
set $S = f$.
\end{PROOF}

\section* {Acknowlegements}
The authors would like to thank D.~Arnal, M.~Cahen, C.~Emmrich, S.~Gutt,
B.~Kostant, and A.~Weinstein for useful discussions.


\begin{thebibliography}{99}

\bibitem {AW70}
         {\sc Agarwal, G. S., Wolf, E.:}
         {\it Calculus for Functions of Noncommuting Operators and General
          Phase-Space Methods in Quantum Mechanics.
          I.~Mapping Theorems and Ordering of Functions of
          Noncommuting Operators.}
          Phys. Rev. D {\bf 2} 10 (1970), 2161-2188.

\bibitem {BFFLS78}
         {\sc Bayen, F., Flato, M., Fronsdal, C.,
         Lichnerowicz, A., Sternheimer, D.:}
         {\it Deformation Theory and Quantization.}
         Annals of Physics {\bf 111}, part I: 61-110,
         part II: 111-151 (1978).

\bibitem {BCG96}
         {\sc Bertelson, M., Cahen, M., Gutt, S.:}
         {\it Equivalence of Star Products.}
         Universit\'e Libre de Bruxelles,
         Travaux de Math\'ematiques, Fascicule {\bf 1}, 1--15
         (1996).

\bibitem {BNW97b}
         {\sc Bordemann, M., Neumaier, N., Waldmann, S.:}
         {\it Homogeneous Fedosov Star Products on Cotangent Bundles II:
         GNS Representations, the WKB Approximation, and Applications.}
         Preprint Univ. Freiburg in preparation.

\bibitem {BW96a}
         {\sc Bordemann, M., Waldmann, S.:}
         {\it A Fedosov Star Product of Wick Type for K\"ahler
         mani\-folds.}
         Preprint Univ. Freiburg FR-THEP-96/9, May 1996, and
         q-alq/9605012. To appear in Lett. Math. Phys.

\bibitem {BW96b}
         {\sc Bordemann, M., Waldmann, S.:}
         {\it Formal GNS Construction and States in Deformation
         Quantization.}
         Preprint Univ. Freiburg FR-THEP-96/12, July 1996,
         and q-alg/9607019.

\bibitem {BW96c}
         {\sc Bordemann, M., Waldmann, S.:}
         {\it Formal GNS Construction and WKB Expansion in
         Deformation Quantization.}
         Contribution to the Ascona Meeting on Deformation Theory,
         Symplectic Geometry, and Applications, June 16--22, 1996.
         q-alg/9611004.

\bibitem {CG82}
         {\sc Cahen, M., Gutt, S.:}
         {\it Regular $*$ Representations of Lie Algebras.}
         Lett. Math. Phys. {\bf 6}, 395--404 (1982).

\bibitem {CFS92}
         {\sc Connes, A., Flato, M., Sternheimer, D.:}
         {\it Closed Star Products and Cyclic Cohomology}
         Lett. Math. Phys. {\bf 24}, 1-12 (1992).

\bibitem {DL83a}
         {\sc DeWilde, M., Lecomte, P. B. A.:}
         {\it Star Products on Cotangent Bundles.}
         Lett. Math. Phys. {\bf 7}, 235--241 (1983).

\bibitem {DL83}
         {\sc DeWilde, M., Lecomte, P. B. A.:}
         {\it Existence of star products and of formal deformations
         of the Poisson Lie Algebra of arbitrary symplectic manifolds.}
         Lett. Math. Phys. {\bf 7}, 487-496 (1983).

\bibitem {DL88}
         {\sc DeWilde, M., Lecomte, P. B. A.:}
         {\it Formal Deformations of the Poisson Lie Algebra of a
         Symplectic Manifold and Star Products. Existence,
         Equivalence, Derivations.}
         in:
         {\sc Hazewinkel, M., Gerstenhaber, M. (eds):}
         {\it Deformation Theory of Algebras and Structures
         and Applications.}
         Kluwer, Dordrecht 1988.

\bibitem {Dri83}
         {\sc Drinfel'd, V. G.:}
         {\it On constant quasiclassical solutions of the Yang-Baxter
         quantum equation.} Soviet Math. Dokl. {\bf 28}, 667--671
         (1983).

\bibitem {EW93}
         {\sc Emmrich, C., Weinstein, A.:}
         {\it The differential geometry of Fedosov's
         quantization.}
         Preprint Univ. California, Berkley, Nov. 1993, and
         hep-th/9311094.

\bibitem {Fed94}
         {\sc Fedosov, B.:}
         {\it A Simple Geometrical Construction of
         Deformation Quantization.}
         J. Diff. Geom. {\bf 40}, 213-238 (1994).

\bibitem {Fed96}
         {\sc Fedosov, B.:}
         {\it Deformation Quantization and Index Theory.}
         Akademie Verlag, Berlin 1996.

\bibitem {Gut83}
         {\sc Gutt, S.:}
         {\it An explicit $*$-Product on the cotangent bundle of a Lie
         Group.}
         Lett. Math. Phys. {\bf 7}, 249--258 (1983).

\bibitem {NT95a}
         {\sc Nest, R., Tsygan, B.:}
         {\it Algebraic Index Theorem.}
         Commun. Math. Phys. {\bf 172}, 223--262 (1995).

\bibitem {NT95b}
         {\sc Nest, R., Tsygan, B.:}
         {\it Algebraic Index Theorem for Families.}
         Adv. Math. {\bf 113}, 151--205 (1995).

\bibitem {Pfl95}
         {\sc Pflaum, M. J.:}
         {\it Local Analysis of Deformation Quantization.}
         Ph.D. thesis, Fakult\"at f\"ur Mathematik der
         Ludwig-Maximilians-Universit\"at, M\"unchen, 1995.

\bibitem {Pfl96}
         {\sc Pflaum, M. J.:}
         {\it The Normal Symbol on Riemannian Manifolds.}
         Preprint dg-ga/9612011.

\bibitem {Und78}
         {\sc Underhill, J.:}
         {\it Quantization on a manifold with connection.}
         J. Math. Phys. {\bf 19}, 1932--1935 (1978).

\bibitem {Wid78}
         {\sc Widom, H.:}
         {\it Families of Pseudodifferential Operators.}
         in:
         {\sc Gohberg, I., Kac, M. (eds):}
         {\it Topics in Functional Analysis.}
         Academic Press, New York 1978.

\bibitem {Wid80}
         {\sc Widom, H.:}
         {\it A Complete Symbol Calculus for Pseudodifferential Operators.}
         Bull. Sci. Math. {\bf 104}, 19--63 (1980).

\bibitem {Woo80}
         {\sc Woodhouse, N.:}
         {\it Geometric Quantization.}
         Clarendon Press, Oxford 1980.

\end {thebibliography}
\end {document}